\documentclass[10pt,oneside,a4paper]{amsart}
\usepackage{amsmath,amsthm, latexsym,verbatim}
\usepackage{amsxtra}
\usepackage{pxfonts}
\usepackage[psamsfonts]{amssymb}
\usepackage[applemac]{inputenc}
\usepackage{fancyhdr}
\usepackage{mathrsfs}
\usepackage[bookmarks,colorlinks]{hyperref}

\fancypagestyle{plain}{
	\fancyhead{}
	
}

\newtheorem{definition}{Definition}[section]
\newtheorem{theorem}[definition]{Theorem}
\newtheorem{proposition}[definition]{Proposition}
\newtheorem{lemma}[definition]{Lemma}

\newtheorem{remark}[definition]{Remark}

\def\sref#1{(\ref{#1})}

\newcommand\dcut{d\llap {\raisebox{.9ex}{$\scriptstyle-\!$}}}

\def\Op#1{ {\operatorname{Op}} \left( #1 \right) }
\def\supp{ {\operatorname{supp}} }
\newcommand{\cl}{\textnormal{cl}}
\newcommand{\m}{\mathbf{m}}
\def\Symp#1{ {\mathrm{Sym_{p}}} \left( #1 \right) }

\def\sref#1{(\ref{#1})}
\def\QED{\fbox{\rule{0.5mm}{0mm}\rule{0mm}{0.5mm}}}
\def\sfrac#1#2{\displaystyle\frac{#1}{#2}}
\newenvironment{Pf}[1]{\textit{Proof{#1}}}{\QED}\def\sref#1{(\ref{#1})}
\def\R {\mathbb{R}}       %%% Real numbers
\def\N {\mathbb{N}}       %%% Natural
\def\C {\mathbb{C}}       %%% Complex
\def\Z {\mathbb{Z}}
\def\Rn {\mathbb{R}^{n}}
\def\SG {S}

\def\ELG {EL}

\def\spec{\operatorname{spec}}

\def\cA{{\mathcal A}}

\def\cD{{\mathcal D}}

\def\cF{{\mathcal F}}
\def\cG{{\mathcal G}}
\def\cJ{{\mathcal J}}
\def\cL{{\mathcal L}}
\def\cN{{\mathcal N}}

\def\cQ{{\mathcal Q}}

\def\cS{{\mathcal S}}
\def\cT{{\mathcal T}}

\def\1{\lambda}
\def\2{\Sigma_{2}}

\def\phasef{\Phi}

\def\<{{\langle}}
\def\>{{\rangle}}
\def\norm#1{{\langle} #1 {\rangle}}

\def\l{{\lambda}}

\def\cl{\mathrm{cl}}

\def\Sn{\mathbb{S}^{n-1}}
\def\SX{ {\mathcal{S}} }

\author{Sandro Coriasco}

\address{Dipartimento di Matematica, Università di Torino, Italy
		\newline\indent
		Institut f\"ur Analysis, Gottfried Wilhelm Leibniz Universit\"at Hannover, Germany}

\email{sandro.coriasco@unito.it}

\author{Lidia Maniccia}

\address{C/O S. Coriasco, Dipartimento di Matematica, Università di Torino, Italy}

\title[On the Spectral Asymptotics of Operators on Manifolds with Ends]{On the Spectral Asymptotics\\of Operators on Manifolds with Ends}



\begin{document}

\begin{abstract}
We deal with the asymptotic behaviour for $\lambda\to+\infty$
of the counting function $N_P(\lambda)$ of certain positive selfadjoint
operators $P$ with double order $(m,\mu)$, $m,\mu>0$, $m\not=\mu$,
defined on a manifold with ends $M$.
The structure of this class of noncompact manifolds allows
to make use of calculi of pseudodifferential operators and
Fourier Integral Operators associated with weighted
symbols globally defined on $\R^n$. By means of these tools,
we improve known results concerning the remainder terms
of the Weyl Formulae for $N_P(\lambda)$ and show how
their behaviour depends on the ratio $\frac{m}{\mu}$ and the
dimension of $M$.
\end{abstract}

\maketitle

\noindent
                         %%%%%%%%%%%%%%%%%%%%%%%%%%%%%%%%%%%
\section{Introduction}
\label{sec:intro}
\setcounter{equation}{0}

The aim of this paper is to study the asymptotic behaviour, for
$\lambda\to+\infty$, of the counting function 
\[
N_P(\lambda) = \sum_{\lambda_j\le\lambda}1
\]
where $\lambda_1\le\lambda_2\le\dots$ is the sequence of the eigenvalues,
repeated according to their multiplicities, 
of a positive order, selfadjoint, classical, elliptic
$SG$-pseudodifferential operator $P$ on a manifold with ends.
Explicitly, $SG$-pseudodifferential operators $P=p(x,D)=\Op{p}$ on
$\R^n$ can be defined via the usual left-quantization
\[
	Pu(x)= \frac{1}{(2 \pi)^{n}} \int e^{i x \cdot \xi} p(x, \xi) \hat u(\xi) d\xi,\quad u\in\SX(\R^n),
\]
starting from symbols $p(x,\xi) \in C^\infty(\R^n\times\R^n)$ with the property that, for arbitrary multiindices $\alpha,\beta$, there exist constants $C_{\alpha\beta}\ge0$ such that the estimates 
\begin{equation}
	\label{disSG}
	|D_\xi^{\alpha}D_x^{\beta} p(x, \xi)| \leq C_{\alpha\beta} \langle\xi\rangle^{m-|\alpha|}\langle x\rangle^{\mu-|\beta|}
\end{equation}
hold for fixed $m,\mu\in\R$ and all $(x, \xi) \in \R^n \times \R^n$, where $\langle y \rangle=\sqrt{1+|y|^2}$, $y\in\R^n$.
Symbols of this type belong to the class denoted by $S^{m,\mu}(\Rn)$, and the corresponding operators constitute the class
$L^{m,\mu}(\Rn)=\Op{S^{m,\mu}(\Rn)}$. In the sequel we will sometimes write $S^{m,\mu}$ and $L^{m,\mu}$, respectively,
fixing once and for all the dimension of the (non-compact) base manifold to $n$. 

These classes of operators, introduced on $\R^n$ by H.O.~Cordes \cite{CO} and C.~Parenti \cite{PA72}, see also R.~Melrose \cite{ME}, M.A.~Shubin \cite{SH87}, form a graded algebra, i.e., $L^{r,\rho}\circ L^{m,\mu}\subseteq L^{r+m,\rho+\mu}$. The remainder elements are operators with symbols in 
$\displaystyle S^{-\infty, -\infty}(\R^n)= \bigcap_{(m,\mu) \in \R^2} S^{m,\mu} (\R^n)=\SX(\R^{2n})$, that is, those having kernel in $\SX(\R^{2n})$, continuously mapping $\SX^\prime(\R^n)$ to $\SX(\R^n)$. An operator $P=\Op{p}\in L^{m,\mu}$ and its symbol $p\in\SG^{m,\mu}$ are called $SG$-elliptic if there exists $R\ge0$ such that $p(x,\xi)$ is invertible for $|x|+|\xi|\ge R$ and
\[
	p(x,\xi)^{-1}=O(\norm{\xi}^{-m}\norm{x}^{-\mu}).
\] 
In such case we will usually write $P\in\ELG^{m,\mu}$.
Operators in $L^{m,\mu}$ act continuously from $\SX(\R^n)$ to itself, and extend as continuous operators from $\SX^\prime(\R^n)$ to itself and from $H^{s,\sigma}(\R^n)$ to $H^{s-m,\sigma-\mu}(\R^n)$, where $H^{s,\sigma}(\R^n)$, $s,\sigma\in\R$, denotes the weighted Sobolev space
\begin{align*}
  	H^{s,\sigma}(\R^n)&= \{u \in \SX^\prime(\R^{n}) \colon \|u\|_{s,\sigma}= \|\Op{\pi_{s,\sigma}}u\|_{L^2}< \infty\},
 	\\
  	\pi_{s,\sigma}(x,\xi)&= \langle \xi \rangle^{s} \langle x\rangle^{\sigma}.
\end{align*}
Continuous inclusions $H^{s,\sigma}(\R^n)\hookrightarrow H^{r,\rho}(\R^n)$ hold when $s\ge r$ and $\sigma\ge\tau$, compact when both inequalities are strict, and
\[
	\displaystyle\SX(\R^n)=\bigcap_{(s,\sigma) \in \R^2} H^{s,\sigma}(\R^n)
	,\qquad
	\displaystyle\SX^\prime(\R^n)=\bigcup_{(s,\sigma) \in \R^2} H^{s,\sigma}(\R^n).
\]
An elliptic $SG$-operator $P \in L^{m,\mu}$ admits a parametrix $E\in L^{-m,-\mu}$ such that
\[
PE=I + K_1, \quad EP= I+ K_2,
\]
for suitable $K_1, K_2 \in L^{-\infty,-\infty}=\Op{S^{-\infty,-\infty}}$, and it turns out to be a Fredholm operator.
In 1988, E.~Schrohe \cite{SC87} introduced a class of non-compact manifolds, the so-called $SG$-manifolds, on which it is possible to transfer from $\R^n$ the whole $SG$-calculus. In short, these are manifolds which admit a finite atlas whose changes of coordinates behave like symbols of order $(0,1)$ (see \cite{SC87} for details and additional technical hypotheses). The manifolds with cylindrical ends are a special case of $SG$-manifolds, on which also the concept of $SG$-classical operator makes sense: moreover, the principal symbol of a $SG$-classical operator $P$ on a manifold with cylindrical ends $M$, in this case a triple $\sigma(P)=(\sigma_\psi(P),\sigma_e(P),\sigma_{\psi e}(P))=(p_\psi, p_e, p_{\psi e})$, has an invariant meaning on $M$, see Y.~Egorov and B.-W.~Schulze \cite{ES97}, L.~Maniccia and P.~Panarese \cite{MP02}, R.~Melrose \cite{ME} and Section \ref{sec:mwe} below. We indicate the subspaces of classical symbols and operators adding the subscript $_\cl$ to the notation introduced above.

The literature concerning the study of the eigenvalue asymptotics of elliptic operators is vast, 
and covers a number of different situations
and operator classes, see, e.g., the monograph by V.J.~Ivrii \cite{Iv98}. 
Then, we only mention a few of the many existing papers
and books on this deeply investigated subject, which are related to the case
we consider here, either by the type of symbols and underlying spaces,
or by the techniques which are used: we
refer the reader to the corresponding reference lists for more
complete informations. On compact manifolds, well known results were proved by 
L.~H\"ormander \cite{Ho2} and V.~Guillemin \cite{GU85},
see also the book by H.~Kumano-go \cite{Kumano-go:1}.
On the other hand, for operators globally defined on $\R^n$, see 
P.~Boggiatto, E.~Buzano, L.~Rodino \cite{BBR96},
B.~Helffer \cite{Helffer:984.1}, 
L.~H{\"o}rmander \cite{Ho3}, A.~Mohammed \cite{Mo1}, F.~Nicola
\cite{NI03}, M.~A.~Shubin \cite{SH87}.
Many other situations have been considered, see the cited book by V.J.~Ivrii.
On manifolds with ends, T.~Christiansen and M.~Zworski \cite{ChZw95}
studied the Laplace-Beltrami operator associated with
a scattering metric, while L.~Maniccia and P.~Panarese \cite{MP02} applied the
heat kernel method to study operators similar to those considered here.

Here we deal with the case of manifolds with ends for $P\in EL^{m,\mu}_\cl(M)$, positive and selfadjoint, such that $m,\mu>0$, $m\not=\mu$,
focusing on the (invariant) meaning of the constants appearing in the 
corresponding Weyl formulae and on achieving a better estimate of the remainder term.
Note that the situation we consider here is different from that of the Laplace-Beltrami
operator investigated in \cite{ChZw95}, where continuous spectrum is present as well:
in fact, in view of Theorem \ref{thm:spt}, $\spec(P)$ consists only of a sequence of real isolated eigenvalues $\{\l_j\}$ with finite multiplicity.

As recalled above, a first result concerning the asymptotic behaviour of 
$N_P(\lambda)$ for operators including those considered in this paper
was proved in \cite{MP02}, giving, for
$\lambda \to +\infty$,
$$
N_P(\lambda)=
\left\{
\begin{array}{ll}
C_1 \lambda^{\frac{n}{m}}+o(\lambda^\frac{n}{m})& \text{for } m<\mu
\\
C_0^1 \lambda^{\frac{n}{m}}\log\lambda+o(\lambda^\frac{n}{m}\log\lambda)& \text{for } m=\mu
\\
C_2 \lambda^{\frac{n}{\mu}}+o(\lambda^\frac{n}{\mu})& \text{for } m>\mu.
\end{array}
\right.
$$
Note that the constants $C_1$, $C_2$, $C_0^1$ above depend only on the principal symbol of $P$, which implies that they have an invariant meaning on the manifold $M$, see Sections \ref{sec:mwe} and \ref{sec:N_a} below.
On the other hand, in view of the technique used there, the remainder
terms appeared in the form $o(\lambda^{\frac{n}{\min\{m,\mu\}}})$ and
$o(\lambda^{\frac{n}{m}}\log\lambda)$ for $m\not=\mu$ and $m=\mu$, respectively.
An improvement in this direction for operators on $\R^n$ had been achieved by F. Nicola
\cite{NI03}, who, in the case $m=\mu$ proved that
\[
N_P(\lambda)=C_0^1\lambda^\frac{n}{m}\log\lambda+O(\lambda^\frac{n}{m}),
\lambda\to+\infty,
\]
while, for $m\not=\mu$, showed that the remainder term has the form
$O(\lambda^{\frac{n}{\min\{m,\mu\}}-\varepsilon})$ for a suitable $\varepsilon>0$.
A further improvement of these results in the case
$m=\mu$ has recently appeared in U. Battisti and S. Coriasco \cite{BC10}, where 
it has been shown that, for a suitable $\varepsilon>0$, 
\[
N_P(\lambda)=C^1_0 \lambda^{\frac{n}{m}}\log \lambda + 
C^2_0 \lambda^{\frac {n}{m}}+ O(\lambda^{\frac{n}{m}-\varepsilon}),
\lambda\to+\infty.
\]
Even the constant $C_0^2$ has an invariant meaning on $M$, and both $C_0^1$ and $C_0^2$ are explicitly computed in terms of trace operators defined on $L^{m,m}_\cl(M)$.

In this paper the remainder estimates in the case $m\neq\mu$ are further improved. 
More precisely, we first consider the power $Q=P^\frac{1}{\max\{m,\mu\}}$ of $P$ (see L. Maniccia, E. Schrohe, J. Seiler \cite{MSS06} for the properties of
powers of $SG$-classical operators). Then, by studying the
asymptotic behaviour in $\lambda$ of the trace of the operator
$\widehat{\psi}_{\lambda}(-Q)$, $\psi_{\lambda}(t)=\psi(t)e^{-it\lambda}$, $\psi \in C_0^\infty(\R)$, defined via a Spectral Theorem
and approximated in terms of Fourier Integral Operators,
we prove the following

\begin{theorem}
\label{thm:main}
Let $M$ be a manifold with ends of dimension $n$ and let $P\in EL^{m,\mu}_\cl(M)$ be a positive selfadjoint operator such that $m,\mu>0$, $m\not=\mu$, with domain $H^{m,\mu}(M)\hookrightarrow L^2(M)$.
Then, the following Weyl formulae hold for $\lambda\to+\infty$:
\begin{equation}
\label{weylnor}
N_P(\lambda) =\left\{
\begin{array}{lll} 
      C_1 \lambda^{\frac{n}{m}}+O(\lambda^\frac{n}{\mu})+O(\lambda^{\frac{n}{m}-\frac{1}{\mu}})
&\!\!\!\!\!= C_1 \lambda^{\frac{n}{m}}+ O(\lambda^{\frac{n}{m}-\varepsilon_1})  & \mbox{for } m<\mu
\\
      C_2 \lambda^{\frac{n}{\mu}}+ O(\lambda^\frac{n}{m})+O(\lambda^{\frac{n}{\mu}-\frac{1}{m}}) 
&\!\!\!\!\!= C_2 \lambda^{\frac{n}{\mu}}+ O(\lambda^{\frac{n}{\mu}-\varepsilon_2}) & \mbox{for }m>\mu.
\end{array}
\right.
\end{equation}
where
$\displaystyle \varepsilon_1=\min\left\{\sfrac{1}{\mu},
n\!\left(\frac{1}{m}-\sfrac{1}{\mu}\right)\right\}$
and 
$\displaystyle \varepsilon_2=\min\left\{\sfrac{1}{m},
n\!\left(\frac{1}{\mu}-\sfrac{1}{m}\right)\right\}$.
\end{theorem}

\noindent
The order of the remainder is then determined by the ratio
of $m$ and $\mu$ and the dimension of $M$, since
\begin{equation}
\label{remorder}
\begin{split}
	\frac{n}{m}-\frac{1}{\mu}\le\frac{n}{\mu} \mbox{ for } m<\mu\Leftrightarrow
	1<\frac{\mu}{m}\le 1+\frac{1}{n},
	 \\
	\frac{n}{\mu}-\frac{1}{m}\le\frac{n}{m} \mbox{ for } m>\mu\Leftrightarrow
	1<\frac{m}{\mu}\le 1+\frac{1}{n}.
\end{split}
\end{equation}
In particular, when $\frac{\max\{m,\mu\}}{\min\{m,\mu\}}\ge2$, the remainder is always $O(\lambda^\frac{n}{\max\{m,\mu\}})$.

Examples include operators of Schr\"odinger type on $M$,
that is $P=-\Delta_g+V$, $\Delta_g$ the Laplace-Beltrami operator in $M$ associated with a suitable metric $g$, $V$ a smooth potential that, in the local coordinates $x\in U_N\subseteq\R^n$ on the cylindrical end growths as $\norm{x}^\mu$, with an appropriate $\mu>0$ related to $g$.
Such examples will be discussed in detail, together with the sharpness of the results in Theorem \ref{thm:main}, in the forthcoming paper \cite{BoC11}, see also \cite{Bo11}.

The key point in the proof of Theorem \ref{thm:main}
is the study of the asymptotic behaviour for $\lambda \to + \infty$ of integrals of the form
$$
I(\l)=\int e^{i(-t \lambda + \varphi(t;x,\xi) - x\xi)}\psi(t) \, a(t;x,\xi)\, dt  \dcut \xi  dx
$$
where $a$ and $\varphi$ satisfy
certain growth conditions in $x$ and $\xi$ (see Section \ref{sec:N_a} for more details). The integrals $I(\lambda)$
represent in fact the local expressions of the trace of
$\widehat{\psi}_{\lambda}(-Q)$, obtained through the so-called ``geometric optic method",
specialised to the SG situation, see e.g. S. Coriasco \cite{Coriasco:998.1, Coriasco:998.2}, S. Coriasco and L. Rodino \cite{Coriasco:998.3}. To treat the
integrals $I(\lambda)$ we proceed similarly to 
A. Grigis and J. Sj\"ostrand \cite{GriSjo:994}, 
B. Helffer and D. Robert \cite{HeRo:1}, see also H. Tamura \cite{Ta82}.

The paper is organised as follows. Section \ref{sec:mwe} is devoted to 
recall the definition of $SG$-classical operators on a manifold with ends $M$. In Section \ref{sec:N_a} we show that the asymptotic behaviour 
of $N_P(\lambda)$, $\lambda\to+\infty$, for a positive self-adjoint operator $P\in L_\cl^{m,\mu}(M)$, $m,\mu>0$, is related to the
asymptotic behaviour of oscillatory integrals of the form $I(\lambda)$. In Section \ref{sec:stat} we conclude the proof of Theorem \ref{thm:main}, investigating the behaviour of $I(\lambda)$ for $\lambda\to+\infty$.
Finally, some technical details are collected in the Appendix.

\section*{Acknowledgements}
The authors wish to thank U. Battisti, L. Rodino and E. Schrohe for useful discussions and hints. Thanks are also due to N. Batavia. The first author was partially supported by the PRIN Project ``Operatori Pseudo-Differenziali ed Analisi Tempo-Frequenza'' (Director of the national project: G. Zampieri; local supervisor at Università di Torino: L. Rodino). The first author also gratefully acknowledges the support by the
Institut für Analysis, Fakult\"at f\"ur Mathematik und Physik, Gottfried Wilhelm Leibniz Universität Hannover,
during his stay as Visiting Scientist in the Academic Year 2011/2012,
where this paper was partly developed and completed.

\section{$SG$-classical operators on manifolds with ends}
\label{sec:mwe}
\setcounter{equation}{0}

From now on, we will be concerned with the subclass of $SG$-operators given by those elements $P\in L^{m,\mu}(\R^n)$, 
$(m,\mu)\in\R^2$, which are $SG$-classical, that is, $P=\Op{p}$ with $p\in \SG^{m,\mu}_\cl(\R^n)\subset \SG^{m,\mu}(\R^n)$. We begin recalling the basic definitions and results (see, e.g., \cite{ES97,MSS06} for additional details and proofs).
\setcounter{equation}{0}

\begin{definition}
\label{def:sgclass-a}
\begin{itemize}
\item[i)]A symbol $p(x, \xi)$ belongs to the class $\SG^{m,\mu}_{\cl(\xi)}(\R^n)$ if there exist $p_{m-i, \cdot} (x, \xi)\in \widetilde{\mathscr{H}}_\xi^{m-i}(\R^n)$, $i=0,1,\dots$, positively homogeneous functions of order $m-i$ with respect to the variable $\xi$, smooth with respect to the variable $x$, such that, for a $0$-excision function $\omega$,
\[
p(x, \xi) - \sum_{i=0}^{N-1}\omega(\xi) \, p_{m-i, \cdot} (x, \xi)\in \SG^{m-N, \mu}(\R^n), \quad N=1,2, \ldots;
\]
\item[ii)]A symbol $p(x, \xi) $ belongs to the class $\SG_{\cl(x)}^{m,\mu}(\R^n)$ if there exist $p_{\cdot, \mu-k}(x, \xi)\in \widetilde{\mathscr{H}}_x^{\mu-k}(\R^n)$, $k=0,\,\dots$, positively homogeneous functions of order $\mu-k$ with respect to the variable $x$, smooth with respect to the variable $\xi$, such that, for a $0$-excision function $\omega$,
\[
p(x, \xi)- \sum_{k=0}^{N-1}\omega(x) \, p_{\cdot, \mu-k}(x,\xi) \in \SG^{m, \mu-N}(\R^n), \quad N=1,2, \ldots
\]
\end{itemize}
\end{definition}
\begin{definition}
\label{def:sgclass-b}
A symbol $p(x,\xi)$ is $SG$-classical, and we write $p \in \SG_{\cl(x,\xi)}^{m,\mu}(\R^n)=\SG_{\cl}^{m,\mu}(\R^n)=\SG_{\cl}^{m,\mu}$, if
\begin{itemize}
\item[i)] there exist $p_{m-j, \cdot} (x, \xi)\in \widetilde{\mathscr{H}}_\xi^{m-j}(\R^n)$ such that, 
for a $0$-excision function $\omega$, $\omega(\xi) \, p_{m-j, \cdot} (x, \xi)\in \SG_{\cl(x)}^{m-j, \mu}(\R^n)$ and
\[
p(x, \xi)- \sum_{j=0}^{N-1} \omega(\xi) \, p_{m-j, \cdot}(x, \xi) \in \SG^{m-N, \mu}(\R^n), \quad N=1,2,\dots;
\]
\item[ii)] there exist $p_{\cdot, \mu-k}(x, \xi)\in \widetilde{\mathscr{H}}_x^{\mu-k}(\R^n)$ such that, 
for a $0$-excision function $\omega$, $\omega(x)\,p_{\cdot, \mu-k}(x, \xi)\in \SG_{\cl(\xi)}^{m, \mu-k}(\R^n)$ and
\[
p(x, \xi) - \sum_{k=0}^{N-1} \omega(x) \, p_{\cdot, \mu-k} \in \SG^{m, \mu-N}(\R^n), \quad N=1,2,\dots
\] 
\end{itemize}
We set $L_{\cl(x, \xi)}^{m,\mu}(\R^n)=L_{\cl}^{m,\mu}=\Op{\SG^{m,\mu}_{\cl}}$.
\end{definition}

\begin{remark}
	The definition could be extended in a natural way
	from operators acting between scalars to operators acting between (distributional sections of) 
	vector bundles: one should then use matrix-valued symbols whose entries satisfy the estimates \eqref{disSG}.
\end{remark}

\noindent
Note that the definition of $SG$-classical symbol implies a condition of compatibility for the terms of the expansions with respect to $x$ and $\xi$. In fact, defining $\sigma_\psi^{m-j}$ and $\sigma_e^{\mu-i}$ on $\SG_{\cl(\xi)}^{m,\mu}$ and $\SG_{\cl(x)}^{m,\mu}$, respectively,  as
\begin{align*}
	\sigma_\psi^{m-j}(p)(x, \xi) &= p_{m-j, \cdot}(x, \xi),\quad j=0, 1, \ldots, 
	\\
	\sigma_e^{\mu-i}(p)(x, \xi) &= p_{\cdot, \mu-i}(x, \xi),\quad i=0, 1, \ldots,
\end{align*}
it possibile to prove that
\[
\begin{split}
p_{m-j,\mu-i}=\sigma_{\psi e}^{m-j,\mu-i}(p)=\sigma_\psi^{m-j}(\sigma_e^{\mu-i}(p))= \sigma_e^{\mu-i}(\sigma_\psi^{m-j}(p)), \\
j=0,1, \ldots, \; i=0,1, \ldots
\end{split}
\]
Moreover, the composition of two $SG$-classical operators is still classical. 
For $P=\Op{p}\in L^{m,\mu}_\cl$ the triple $\sigma(P)=(\sigma_\psi(P),\sigma_e(P),\sigma_{\psi e}(P))=
(p_{m,\cdot}\,,\,p_{\cdot,\mu}\,,\, p_{m,\mu})=(p_\psi,p_e,p_{\psi e})$
is called the \textit{principal symbol of $P$}. The three components are also called the $\psi$-, $e$- and
$\psi e$-principal symbol, respectively. This definition keeps the usual multiplicative behaviour, that is, 
for any $R\in L^{r,\rho}_\cl$, $S\in L^{s,\sigma}_\cl$, $(r,\rho),(s,\sigma)\in\R^2$,
$\sigma(RS)=\sigma(S)\,\sigma(T)$, with componentwise product in the right-hand side. We also set 
\begin{align*}
	\Symp{P}(x,\xi) = &\;\,\Symp{p}(x,\xi) =\\
	=
	&\;\, p_\m(x,\xi)=\omega(\xi) p_{\psi}(x,\xi) +
	\omega(x)(p_{e}(x,\xi) - \omega(\xi) p_{\psi e}(x,\xi)),
\end{align*}
for a fixed $0$-excision function $\omega$. Theorem \ref{thm:ellclass} below allows to express the ellipticity
of $SG$-classical operators in terms of their principal symbol:
\begin{theorem}
	 \label{thm:ellclass}
	An operator $P\in L^{m,\mu}_\cl$ is elliptic if and only if each element of the triple $\sigma(P)$ is 
	non-vanishing on its domain of definition.
\end{theorem}

\noindent
As a consequence, denoting by $\{\lambda_j\}$ the sequence of eigenvalues of $P$, ordered such that $j\le k\Rightarrow \lambda_j\le\lambda_k$, with each eigenvalue repeated accordingly to its multiplicity, the counting function $\displaystyle N_P(\lambda)=\sum_{\lambda_j\le \lambda} 1$ is well-defined for a $SG$-classical elliptic self-adjoint operator $P$,
see, e.g., \cite{BC10, Bo11, BoC11, NI03}. 
%
% Manifolds with ends
%
We now introduce the class of noncompact manifolds with which we will deal:

\begin{definition}
	A manifold with a cylindrical end is a triple $(M, X, [f])$, where $M= \mathscr{M} \amalg_C \mathscr{C}$ is a $n$-dimensional
	smooth manifold and
	\begin{enumerate}
		\item[  i)] $\mathscr{M}$ is a smooth manifold, given by $\mathscr{M}=(M_0\setminus D)\cup C$ 
		with a $n$-dimensional smooth compact manifold without boundary $M_0$, $D$ a closed disc of $M_0$ and
		$C\subset D$ a collar neighbourhood of $\partial D$ in $M_0$;
		\item[ ii)] $\mathscr{C}$ is a smooth manifold with boundary $\partial\mathscr{C}=X$, with $X$ diffeomorphic to 
		$\partial D$;
		\item[iii)] $f: [\delta_f, \infty) \times \mathbb{S}^{n-1} \rightarrow \mathscr{C}$, $\delta_f>0$,
				is a diffeomorphism, $f(\{\delta_f \}\times \mathbb{S}^{n-1})=X$ and 
				$f(\{[\delta_f,\delta_f+\varepsilon_f) \}\times \mathbb{S}^{n-1})$, $\varepsilon_f>0$,
				is diffeomorphic to $C$;  
		\item[ iv)] the symbol $\amalg_C$ means that we are gluing $\mathscr{M}$ and $\mathscr{C}$,
				through the identification of $C$ and $f(\{[\delta_f,\delta_f+\varepsilon_f) \}\times \mathbb{S}^{n-1})$;
		\item[  v)] the symbol $[f]$ represents an equivalence class in the set of functions
				\begin{align*}
					\{ g: [\delta_g, \infty) \times \mathbb{S}^{n-1} \rightarrow \mathscr{C} \colon & g \textnormal{ is a diffeomorphism, } \\
						&g(\{\delta_g\}\times \mathbb{S}^{n-1})=X\mbox{ and } \\
						&g([\delta_g, \delta_g+\varepsilon_g) \times \mathbb{S}^{n-1}),
						\mbox{ $\varepsilon_g>0$, is diffeomorphic to $C$}\}
				\end{align*}
				where $f \sim g$ if and only if there exists a diffeomorphism 
				$\Theta \in \textnormal{Diff}(\mathbb{S}^{n-1})$ such that
				\begin{equation}
					\label{econd}
					(g^{-1} \circ f)(\rho, \gamma)= (\rho, \Theta(\gamma))
				\end{equation}
				for all $\rho\ge \max\{\delta_f, \delta_g\}$ and $\gamma \in \mathbb{S}^{n-1}$. 
	\end{enumerate}
\end{definition}

\noindent
We use the following notation: 
\begin{itemize}
\item $U_{\delta_f}= \{x \in \R^n \colon |x|> \delta_f\}$;
\item $ \mathscr{C}_\tau= f([\tau, \infty) \times \mathbb{S}^{n-1})$, where $\tau\ge\delta_f$.
The equivalence condition \eqref{econd} implies that $\mathscr{C}_\tau$ is well defined;
\item $\displaystyle \pi: \R^n\setminus\{0\}\rightarrow  (0, \infty) \times \mathbb{S}^{n-1}: x \mapsto \pi(x)= \Big(|x|, \frac{x}{|x|}\Big)$;
\item $f_\pi= f\circ \pi: \overline{U_{\delta_f}} \rightarrow \mathscr{C}$ is a parametrisation of the end. Let us notice that, setting
$F=g^{-1}_\pi \circ f_\pi$, the equivalence condition \eqref{econd} implies
\begin{equation}
	\label{cambcart}
	F(x)= |x| \; \Theta\Big(\frac{x}{|x|}\Big).
\end{equation}
We also denote the restriction of $f_\pi$ mapping $U_{\delta_f}$ onto $\dot{\mathscr{C}}=\mathscr{C}\setminus X$ by
$\dot{f}_\pi$.
\end{itemize}
The couple $(\dot{\mathscr{C}}, \dot{f}_\pi^{-1})$ is called the exit chart.
If $\mathscr{A}=\{(\Omega_i, \psi_i)\}_{i=1}^N$ is such that the subset $\{(\Omega_i, \psi_i)\}_{i=1}^{N-1}$ is a finite atlas for 
$\mathscr{M}$ and $(\Omega_N, \psi_N)=(\dot{\mathscr{C}}, \dot{f}_\pi^{-1})$, then $M$, with the atlas $\mathscr{A}$, is a $SG$-manifold (see \cite{SH87}): an atlas $\mathscr{A}$ of such kind is called \textit{admissible}. From now on, we restrict the choice of atlases on $M$ to the class of admissible ones. We introduce the following spaces, endowed with their natural topologies:

\[
\begin{split}
\mathscr{S}(U_\delta)&=\left\{u \in C^{\infty}(U_\delta) \colon \forall \alpha, \beta \in \N^n\, \forall \delta'>\delta \,\sup_{x\in U_{\delta^\prime}}|x^\alpha\partial^\beta u(x)|< \infty\right\},\\
\mathscr{S}_0(U_\delta)&= \bigcap_{\delta' \searrow \delta}\{u \in \mathscr{S}(\R^n)\colon \mathrm{supp}\,u \subseteq \overline{U_{\delta'}} \}, \\
\mathscr{S}(M)&= \{ u \in C^{\infty}(M) \colon  u \circ \dot{f}_\pi \in \mathscr{S}(U_{\delta_f}) \mbox{ for any exit map } f_\pi \},\\
\mathscr{S}^\prime(M)&\mbox{ denotes the dual space of $\mathscr{S}(M)$}.
\end{split}
\]

\begin{definition}
The set $\SG^{m, \mu}(U_{\delta_f})$ consists of all the symbols $a \in C^{\infty}(U_{\delta_f})$ which fulfill \eqref{disSG} for $(x,\xi) \in U_{\delta_f} \times \R^n$ only. Moreover, the symbol $a$ belongs to the subset $SG_{\cl}^{m, \mu}(U_{\delta_f})$ if it admits expansions in asymptotic sums of homogeneous symbols with respect to $x$ and $\xi$ as in Definitions \ref{def:sgclass-a} and \ref{def:sgclass-b}, where the remainders are now given by $SG$-symbols of the required order on $U_{\delta_f}$.
\end{definition}

\noindent
Note that, since $U_{\delta_f}$ is conical, the definition of homogeneous and classical symbol on $U_{\delta_f}$ makes sense. Moreover,
the elements of the asymptotic expansions of the classical symbols can be extended by homogeneity to smooth functions on
$\R^n\setminus\{0\}$, which will be denoted by the same symbols.
It is a fact that, given an admissible atlas $\{(\Omega_i, \psi_i)\}_{i=1}^N$ on $M$, there exists a partition of unity $\{\theta_i\}$ and a
set of smooth functions $\{\chi_i\}$ which are compatible with the $SG$-structure of $M$, that is:
\begin{itemize}
\item $\mathrm{supp}\,\theta_i\subset\Omega_i$, $\mathrm{supp}\,\chi_i\subset\Omega_i$, $\chi_i\,\theta_i=\theta_i$, $i=1,\dots,N$;
\item $|\partial^\alpha(\theta_N\circ \dot{f}_\pi)(x)|\le C_\alpha \norm{x}^{-|\alpha|}$ and 
$|\partial^\alpha(\chi_N\circ \dot{f}_\pi)(x)|\le C_\alpha \norm{x}^{-|\alpha|}$ for all $x\in U_{\delta_f}$.
\end{itemize}
Moreover, $\theta_N$ and $\chi_N$ can be chosen so that $\theta_N\circ \dot{f}_\pi$ and $\chi_N\circ \dot{f}_\pi$
are homogeneous of degree $0$ on $U_\delta$. We denote by
$u^*$ the composition of $u\colon \psi_i(\Omega_i)\subset\R^n\to\C$ with the coordinate patches $\psi_i$,
and by $v_*$ the composition of $v\colon \Omega_i\subset M\to \C$ with $\psi_i^{-1}$, $i=1,\dots,N$.
It is now possible to give the definition of $SG$-pseudodifferential operator on $M$:

\begin{definition}
\label{def:mwce}
Let $M$ be a manifold with a cylindrical end. A linear operator $P:\mathscr{S}(M)\to \mathscr{S}'(M) $ is a $SG$-pseudodifferential operator of order $(m, \mu)$ on $M$,
and we write $P\in L^{m,\mu}(M)$, if, for any admissible atlas $\{(\Omega_i, \psi_i)\}_{i=1}^N$ on $M$ with exit chart
$(\Omega_N,\psi_N)$:
\begin{itemize}
\item[1)] for all $i=1, \ldots, N-1$ and any $\theta_i,\chi_i\in C_c^{\infty}(\Omega_i)$, 
there exist symbols $p^i(x,\xi) \in S^{m}(\psi_i(\Omega_i))$ such that
\[
(\chi_i P \theta_i \,u^*)_*(x)= \iint e^{i(x-y)\cdot \xi}p^i(x,\xi) u(y) dy dx, \quad u \in C^{\infty}(\psi_i(\Omega_i));
\]
\item[2)] for any $\theta_N,\chi_N$ of the type described above, there exists a symbol $p^N(x,\xi) \in SG^{m, \mu}(U_{\delta_f})$ such that
\[
(\chi_N P \theta_N\,u^*)_*(x)= \iint e^{i(x-y)\cdot \xi}p^N(x,\xi) u(y) dy dx, \quad u \in \mathscr{S}_0(U_{\delta_f});
\]
\item[3)] $K_P$, the Schwartz kernel of $P$, is such that
\[
K_P \in C^{\infty}\big((M \times M) \setminus \Delta\big) \bigcap \mathscr{S}\big((\dot{\mathscr{C}} \times \dot{\mathscr{C}})\setminus W\big)
\]
where $\Delta$ is the diagonal of $M \times M$ and $W= (\dot{f}_\pi \times \dot{f}_\pi)(V)$ with any conical neighbourhood $V$ of the diagonal of $U_{\delta_f} \times U_{\delta_f}$.
\end{itemize}
\end{definition}

\noindent
The most important local symbol of $P$ is $p^N$. Our definition of $SG$-classical operator on
$M$ differs slightly from the one in \cite{MP02}:

\begin{definition}
\label{clexit}
Let $P \in L^{m, \mu}(M)$. $P$ is a $SG$-classical operator on $M$, and we write $P \in L_{\cl}^{m, \mu}(M)$, if
$p^N(x,\xi) \in \SG_{\cl}^{m, \mu}(U_{\delta_f})$ and the operator $P$, restricted to the manifold
$\mathscr{M}$, is classical in the usual sense.
\end{definition}

\noindent
The usual homogeneous principal symbol $p_{\psi}$ of a $SG$-classical operator $P\in L^{m,\mu}_\cl(M)$ is of course well-defined as a smooth function on $T^*M$. In order to give an invariant definition of the principal symbols homogeneous
in $x$ of an operator $P \in L^{m, \mu}_\cl(M)$, the subbundle $T_X^*M=\{(x,\xi) \in T^*M\colon x \in X, \, \xi \in T_x^*M\}$ was introduced. The notions of ellipticity can be extended to operators on $M$ as well:

\begin{definition}
Let $P \in L_\cl^{m, \mu}(M)$ and let us fix an exit map $f_\pi$. We can define
local objects $p_{m-j, \mu-i}, p_{\cdot, \mu-i} $ as
\[
\begin{split}
p_{m-j, \mu-i}(\theta, \xi)&= p^N_{m-j, \mu-i}(\theta, \xi), \quad \theta \in \mathbb{S}^{n-1}, \,\xi \in \R^n \setminus\{0\},\\
p_{\cdot, \mu-i}(\theta,\xi)&=p^N_{\cdot, \mu-i}(\theta, \xi), \quad \theta \in \mathbb{S}^{n-1},\, \xi \in \R^n.
\end{split}
\]
\end{definition}

\begin{definition}
An operator $P \in L^{m, \mu}_\cl(M)$ is elliptic, and we write $P \in EL^{m, \mu}_\cl(M)$, if the principal part of $p^N \in \SG^{m, \mu}(U_{\delta_f})$ satisfies the $SG$-ellipticity conditions on $U_{\delta_f}\times\R^n$ and the operator $P$, restricted to the manifold
$\mathscr{M}$, is elliptic in the usual sense.
\end{definition}

\begin{proposition}
\label{prop:classinv}
The properties $P \in L^{m, \mu}(M)$ and $P\in L^{m, \mu}_\cl(M)$, as well as the notion of $SG$-ellipticity, do not depend on the (admissible) atlas on $M$. Moreover, the local functions $p_{e}$ and $p_{\psi e}$ give rise to invariantly defined elements of $C^{\infty}(T_X^*M)$ and $C^{\infty}(T_X^*M\setminus 0)$, respectively.
\end{proposition}
\noindent
Then, with any $P\in L^{m,\mu}_\cl(M)$, it is associated an invariantly defined principal symbol in three components
$\sigma(P)=(p_{\psi},p_{e},p_{\psi e})$. Finally, through local symbols given by $\pi^j_{s,\sigma}(x,\xi)=\norm{\xi}^{s}$,
$j=1,\dots,N-1$, and $\pi^N_{s,\sigma}(x,\xi)=\norm{\xi}^{s}\norm{x}^{\sigma}$, $s,\sigma\in\R$,  we get a 
$SG$-elliptic operator
$\Pi_{s,\sigma}\in L^{s,\sigma}_\cl(M)$ and introduce the (invariantly defined) weighted Sobolev spaces $H^{s,\sigma}(M)$ as
\[
	H^{s,\sigma}(M)=\{u\in\mathscr{S}^\prime(M)\colon \Pi_{s,\sigma}u\in L^2(M)\}.
\] 
The properties of the spaces $H^{s,\sigma}(\R^n)$ extend to $H^{s,\sigma}(M)$ without any change, as well as the continuity of the linear mappings $P\colon H^{s,\sigma}(M)\to H^{s-m,\sigma-\mu}(M)$ induced by $P\in L^{m,\mu}(M)$, mentioned in Section \ref{sec:intro}.

\section{Spectral  asymptotics for $SG$-classical elliptic self-adjoint operators\\ on manifolds with ends}
\label{sec:N_a}
\setcounter{equation}{0}

In this section we illustrate the procedure to prove Theorem \ref{thm:main}, similarly to \cite{GriSjo:994}, \cite{Helffer:984.1}, \cite{Ta82}. The result will follow from the Trace formula \sref{eq:tr2}, \sref{eq:tr3}, the asymptotic behaviour \sref{eq:asint1} and the Tauberian Theorem \ref{thm:taub}.
The remaining technical points, in particular the proof of the asymptotic behaviour of the integrals appearing in \sref{eq:tr3}, are described in Section \ref{sec:stat} and in the Appendix. 
 
 Let the operator $P \in \ELG^{m,\mu}_\cl(M)$ be considered as an unbounded operator
$P\colon \cS(M) \subset H^{0,0}(M)=L^2(M) \rightarrow L^2(M)$.
The following Proposition can be proved by reducing to the local situation and using continuity and ellipticity of $P$, its parametrix and the density of $\cS(M)$ in
the $H^{s,\sigma}(M)$ spaces,

\begin{proposition}
	\label{prop:spp}
	Every $P \in \ELG^{m,\mu}_\cl(M)$, considered as an unbounded operator
	$P\colon \cS(M) \subset L^2(M) \rightarrow L^2(M)$, admits a unique closed extension, still denoted
	by $P$, whose domain is $\cD(P) = H^{m,\mu}(M)$.
\end{proposition}

\noindent
From now on, when we write $P\in\ELG^{m,\mu}_\cl(M)$ we always mean its unique closed extension, defined in Proposition \ref{prop:spp}. As standard, we denote by $\varrho(P)$ the resolvent set of $P$, i.e., the set of all $\lambda\in\C$ such that $\lambda I - P$ maps $H^{m,\mu}(M)$ bijectively onto $L^2(M)$. The spectrum of $P$ is then $\spec(P)=\C\setminus\varrho(P)$. The next Theorem was proved in \cite{MP02}.

\begin{theorem}
	\label{thm:spt}
	(Spectral Theorem)
Let $P\in \ELG^{m,\mu}_\cl(M)$ be regarded as a closed unbounded operator on $L^{2}(M)$ with dense
domain $H^{m,\mu}(M)$. Assume also that $m,\mu>0$ and $P^*=P$. Then:
    \begin{itemize}
    \item[i)] $(\lambda I-P)^{-1}$ is a compact operator on $L^{2}(M)$ for
        every $\lambda \in \varrho(P)$. More precisely, $(\lambda
        I-P)^{-1}$ is an extension by continuity from $\cS(M)$
        or a restriction from $\cS^\prime(M)$ of an operator in
        $\ELG^{-m,-\mu}_\cl(M)$.
    \item[ii)] $\spec(P)$ consists
        of a sequence of real isolated eigenvalues $\{\l_j\}$ with finite multiplicity,
        clustering at infinity; the orthonormal system
        of eigenfunctions $\{e_j\}_{j\geq 1}$ is complete in $L^2(M)=H^{0,0}(M)$. Moreover,
        $e_j\in \cS(M)$ for all $j$.
    \end{itemize}
\end{theorem}

Given a positive selfadjoint operator $P\in\ELG_\cl^{m,\mu}(M)$,  $m,\mu>0$, 
$\mu\not= m$, we can assume, without loss of generality, $1 \le \lambda_1 \le \lambda_2 \dots$
(considering, if necessary, $P+c$ in place of $P$, with $c\in\R$ a suitably large constant).
Define the counting function $N_P(\lambda)$, $\lambda\in\R$, as
\begin{equation}
	\label{eq:NP}
	N_P(\lambda) =  \sum_{\lambda_j\le\lambda}1 = \#(\spec(P) \cap (-\infty,\lambda]).
\end{equation}
Clearly, $N_P$ is non-decreasing, continuous from the right and supported in $[0,+\infty)$.
If we set $Q=P^\frac{1}{l}$, 
$l=\max\{m,\mu\}$ (see \cite{MSS06} for the definition of the powers of $P$), $Q$ turns out to be a $SG$-classical elliptic selfadjoint operator
with $\sigma(Q)=(p_\psi^\frac{1}{l},
p_e^\frac{1}{l},p_{\psi e}^\frac{1}{l})$. We denote 
by $\{ \eta_j \}$ the sequence of eigenvalues of $Q$, which satisfy 
$\eta_j=\lambda_j^\frac{1}{l}$: we can then, as above, consider $N_Q(\eta)$.
It is a fact that $N_Q(\eta)=O(\eta^\frac{n}{l})$, see \cite{MP02}. 

From now on we focus on the case $\mu>m>0$: the case $m>\mu>0$ can be treated in a completely similar way, exchanging the role of $x$ and $\xi$. So we can start from a closed
positive selfadjoint operator $Q\in\ELG_\cl^{m,1}(M)$ 
with domain $\cD(Q)=H^{m,1}(M)$, $m\in(0,1)$. 
For $u \in H^{m,1}(M)$, $t\in\R$, we set
\begin{equation}
	\label{eq:U}
	U(t) u = \sum_{j=1}^\infty e^{it\eta_j} \, (u, e_j)_{L^2(M)} \, e_j,
\end{equation}
and the series converges in the $L^2(M)$ norm (cfr., e.g., \cite{GriSjo:994}).
Clearly, for all $t\in \R$, $U(t)$ is 
a unitary operator such that
\[
U(0) = I, \hspace{0.5cm} U(t+s)=U(t) \, U(s), \hspace{0.5cm} t,s \in 
\R.
\]
Moreover, if $u\in H^{km,k}(M)$ for some $k\in\N$, $U(t)u \in 
C^k(\R,H^{0,0}(M)) \cap \ldots \cap C^0(\R,H^{km,k}(M))$ and,
for $u\in H^{m,1}(M)$, we have $D_t U(t) u - Q U(t) u = 0$, 
$U(0)u=u$, which implies that $v(t,x)=U(t)\,u(x)$ is a solution of 
the Cauchy problem
\begin{equation}
	\label{eq:CP}
	(D_t - Q) v = 0, \hspace{0.5cm} \left. v \right|_{t=0} = u.
\end{equation}
Let us fix $\psi \in \cS(\R)$. We can then define the operator 
$\widehat{\psi}(-Q)$ either by using the formula
\[
\widehat{\psi}(-Q)u = \sum_{j=1}^\infty \widehat{\psi}(-\eta_j) \, 
(u,e_j)_{L^2(M)} \, e_j
\]
or by means of the vector-valued integral $\displaystyle \left( \int 
\psi(t) \, U(t) dt \right) u =
\int \psi(t)\,U(t)\,u \,dt$, $u\in H^{0,0}(M)$. Indeed, there exists 
$N_0 \in \N$ such that $\displaystyle \sum_{j=1}^\infty \eta_j^{-N_0} 
< \infty$, so the definition makes sense and gives an operator in 
$\cL(L^2(M))$ with norm bounded by $\| \psi \|_{L^1(\R)}$. The 
following Lemma, whose proof can be found in the Appendix,
is an analog on $M$ of Proposition 1.10.11 in 
\cite{Helffer:984.1}:

\begin{lemma}
	\label{lemma:kernconv}
	$\widehat{\psi}(-Q)$ is an operator with kernel $K_\psi(x,y) = \displaystyle \sum_j \widehat{\psi}(-\eta_j)
	e_j(x) \overline{e_j(y)} \in \cS(M\times M)$.
\end{lemma}

\noindent
Clearly, we then have
\begin{equation}
	\label{eq:tr1}
	\int_M K_\psi(x,x) \, dx = \sum_j \widehat{\psi}(-\eta_j).
\end{equation} 

By the analysis in \cite{Coriasco:998.1, Coriasco:998.2}, 
\cite{Coriasco-Panarese:001}, \cite{Coriasco:998.3} (see also \cite{Coriasco-Maniccia:003}), the above Cauchy Problem \eqref{eq:CP} can be solved modulo $\cS(M)$ by means of a smooth family of operators
$V(t)$, defined for $t \in (-T, T )$, 
$T > 0$ suitably small, in the sense that 
$(D_t-Q)\circ V$ is a family of smoothing operators and $V(0)$ is the identity
on $\cS^\prime(M)$.
More explicitly, 
the following theorem holds (see the Appendix for some details concerning the extension to
the manifold $M$ of the results on $\R^n$ proved in \cite{Coriasco:998.1, Coriasco:998.2}, 
\cite{Coriasco-Panarese:001}, \cite{Coriasco:998.3}).
\begin{theorem}
	\label{thm:globalV}
	Define
	$\displaystyle V(t)u = \sum_{k=1}^N \chi_k A_{k}(t)( \theta_k u)$, 
	where $\theta_k$ and $\chi_k$ are as in Definition \ref{def:mwce}, with
	$\chi_k\,\theta_k=\theta_k$, $k=1,\dots,N$, while the $A_k(t)$ are $SG$ 
	FIOs which, in the local coordinate open set $U_k=\psi_k(\Omega_k)$ and with
	$v \in \cS(\R^n)$, are given by
	$$
		(A_k(t)v)(x) = \int e^{i \varphi_k(t; x, \xi)} \, a_k(t; x, \xi) \, 
		{\hat{v}}(\xi) \,\dcut \xi.
	$$
	Each $A_k(t)$ solves a local Cauchy Problem 
	$(D_t - Q_k) \circ A_k \in C^\infty((-T,T), L^{-\infty,-\infty}(\R^n))$, $A_k(0)=I$,
	with  $Q_k=\mathrm{op}\,(q_k)$ and
	$\{ q_k \} \subset SG_\cl^{m,1}(\R^n)$ local (complete) symbol of $Q$ associated with 
	$\{ \theta_k \}$, $\{ \chi_k \}$, with phase and amplitude functions such that
	\begin{equation}
		\label{eq:ph-amp}
		\begin{aligned}
		& \partial_t \varphi_k(t;x,\xi) - q_k(x,d_x \varphi_k(t;x,\xi)) = 0, \varphi_k(0;x,\xi) = x \xi,
		\\
		& a_k \in C^\infty((-T,T), SG_\cl^{0,0}(\R^n)), a_k(0;x,\xi)=1.
		\end{aligned}
	\end{equation}
	Then, $V(t)$ satisfies 
	$$
	(D_t-Q)\circ V \in C^\infty((-T,T), L^{-\infty,-\infty}(M)),
	\hspace{0.5cm}
	V(0) = I,
	$$
	and
	$U - V \in C^\infty((-T,T), L^{-\infty,-\infty}(M))$.
\end{theorem}

\begin{remark}
	Trivially, for $k=1, \dots, N-1$,
	$q_k$ and $a_k$
	can be considered $SG$-classical, since, in those cases, they actually have order $-\infty$
	with respect to $x$, by the fact that $q_k(x,\xi)$ vanishes for $x$ outside a compact set.
\end{remark}

\begin{remark}
	Notation like $b\in C^\infty((-T,T), \SG^{r,\rho}(\R^n))$, $B\in C^\infty((-T,T), L^{r,\rho}(M))$,
	and similar,
	in Theorem \ref{thm:globalV} and in the sequel,
	also mean that the seminorms of the involved elements 
	in the corresponding spaces (induced, in the mentioned cases, by \eqref{disSG}), are uniformly bounded
	with respect to $t\in (-T,T)$.
\end{remark}

If we write $\psi_\lambda(t)=\psi(t) e^{-i t\lambda}$ in place of  $\psi(t)$, for a chosen $\psi \in C_0^\infty((-T,T))$, the trace formula \sref{eq:tr1} becomes
\begin{equation}
	\label{eq:tr2}
	\int_M K_{\psi_\lambda}(x,x) \, dx = \sum \widehat{\psi}(\lambda-\eta_j).
\end{equation}
Let us denote the kernel of $U-V$ by $r(t;x,y)\in 
C^\infty((-T,T), \cS(M\times M))$. Then, the distribution kernel of
$\displaystyle \int e^{-it\lambda} 
\, \psi(t) \, U(t) \, dt = \widehat{\psi}_\lambda(-Q)$ is
\begin{align*}
K_{\psi_\lambda}(x,y) &=
\sum_{k=1}^N \chi_k(x) \int \!\! \int \psi(t) \,
e^{i(-t\lambda + \varphi_k(t;x,\xi) - y \xi)}
a_k(t;x,\xi)  \,dt \, \dcut \xi \, \theta_k(y) 
\\
&+ \int e^{-it\lambda} \, \psi(t) \, r(t;x,y) \, dt,
\end{align*}
where the local coordinates in the right hand side depend on $k$ and, to simplify the notation, we have omitted the corresponding coordinate maps.
By the choices of $\psi$, $\theta_k$ and $\chi_k$ we obtain
\begin{eqnarray}
\label{eq:tr3}
\sum_j \widehat{\psi}(\lambda-\eta_j)
& = & \sum_{k=1}^N \int \!\! \int \!\! \int \psi(t) \,
e^{i(-t\lambda + \varphi_k(t;x,\xi) - x \xi)}
a_k(t;x,\xi) \, \theta_k(x) \, dt  \, \dcut\xi \, dx
\\
\nonumber
& + & \int\!\int e^{-it\lambda} \, \psi(t) \, r(t;x,x)  \, dt \, dx
\\
\nonumber
& = & \sum_{k=1}^N \int \!\! \int \!\! \int \psi(t) \,
e^{i(-t\lambda + \varphi_k(t;x,\xi) - x \xi)}
a_k(t;x,\xi) \, \theta_k(x) \, dt \, \dcut\xi \, dx
\\
\nonumber
& + &O(|\lambda|^{-\infty}).
\end{eqnarray}
Let $\psi \in C_0^\infty((-T,T))$, $T>0$, be such that $\psi(0)=1$ and $\widehat{\psi}\ge 0$, $\widehat{\psi}(0)>0$ (e.g., set $\psi=\chi * \overline{\check{\chi}}$
with a suitable $\chi\in C_0^\infty((-T,T))$). By the analysis of the asymptotic behaviour of the integrals appearing in \sref{eq:tr3}, described in Section \ref{sec:stat}, we finally obtain
\begin{equation}
	\label{eq:asint1}
	\sum_{j} \widehat{\psi}(\lambda-\eta_j) = 
	\left\{
	\begin{array}{ll}
		\displaystyle \sfrac{n}{m}d_0 \, \lambda^{\frac{n}{m}-1} + O(\lambda^{n^*-1}) &
		 \mbox{for } \lambda\to+\infty
		\\
		\displaystyle\rule{0mm}{7mm}O(|\lambda|^{-\infty}) & \mbox{for } \lambda\to-\infty,
	\end{array}
	\right.
\end{equation}
with $n^* =\min\left\{ n, \sfrac{n}{m} - 1 \right\}$.
The following Tauberian Theorem is a slight modification of Theorem 4.2.5
of \cite{Helffer:984.1} (see the Appendix):

\begin{theorem}
	\label{thm:taub}
	Assume that
	\begin{enumerate}
		\item[  i)] $\psi \in C_0^\infty(\R)$ is an even function satisfying
		                $\psi(0)=1$, $\widehat{\psi}\ge 0$, $\widehat{\psi}(0)>0$;
		\item[ ii)] $N_Q(\lambda)$ is a nondecreasing function, supported in $[0,+\infty)$, continuous from the
		                right, with polynomial growth at infinity and isolated discontinuity points of first kind
		                $\{ \eta_j \}$, $j\in\N$, such that $\eta_j\to+\infty$;
		 \item[iii)] there exists $d_0 \ge 0$ such that
		$$
		\sum_{j} \widehat{\psi}(\lambda-\eta_j) =
		\int \widehat{\psi}(\lambda-\eta) dN_Q(\eta) = 
		\left\{
		\begin{array}{ll}
		\displaystyle \sfrac{n}{m}d_0 \, \lambda^{\frac{n}{m}-1} + O(\lambda^{n^*-1}) &
		 \mbox{for } \lambda\to+\infty
		\\
		\displaystyle\rule{0mm}{7mm}O(|\lambda|^{-\infty}) & \mbox{for } \lambda\to-\infty,
		\end{array}
		\right.
		$$
		with $m\in(0,1)$, $n^* =\min\left\{ n, \sfrac{n}{m} - 1 \right\}$. 
	\end{enumerate}
	Then
	$$
	N_Q(\lambda) = \frac{d_0}{2\pi} \lambda^\frac{n}{m} + O(\lambda^{n^*}),
	\mbox{ for } \lambda\to+\infty.
	$$
\end{theorem}
\begin{remark}
	\label{rem:asint}
	The above statement can be modified as follows: with $\psi$, $N_Q$ and $m$ as in 
	Theorem \ref{thm:taub}, when
	$$
		\int \widehat{\psi}(\lambda-\eta) dN_Q(\eta) = 
		\left\{
		\begin{array}{ll}
		\displaystyle \sfrac{n}{m}d_0 \, \lambda^{\frac{n}{m}-1} + O(\lambda^{\frac{n}{m}-2})
		+ O(\lambda^{n-1}) &
		 \mbox{for } \lambda\to+\infty
		\\
		\displaystyle\rule{0mm}{7mm}O(|\lambda|^{-\infty}) & \mbox{for } \lambda\to-\infty,
		\end{array}
		\right.
	$$
	with $m\in(0,1)$, then
	$\displaystyle
	N_Q(\lambda) = \frac{d_0}{2\pi} \lambda^\frac{n}{m} + O(\lambda^{\frac{n}{m}-1})
		+ O(\lambda^{n})$,
	for $\lambda\to+\infty$.
\end{remark}

\section{Proof of Theorem \ref{thm:main}}
\label{sec:stat}
\setcounter{equation}{0}

In view of Theorem \ref{thm:taub} and Remark \ref{rem:asint}, to complete the proof of Theorem \ref{thm:main} we need to show that \eqref{eq:asint1} holds. To this aim,
as explained above, this Section will be devoted to studying the asymptotic behaviour for
$|\l| \to + \infty$ of 
\begin{equation}
\label{I}
I(\l)=\int e^{i\phasef(t;x,\xi;\l)}\psi(t) \, a(t;x,\xi)\, dt  \dcut \xi  dx,
\end{equation}
where $\psi \in C^\infty_0((-T,T))$, $\psi(0)=1$, $a\in  C^\infty((-T,T),S^{0,0}(\R^n))$, $a(0;x,\xi)=1$, and
$$
\phasef(t;x,\xi;\l)=\varphi(t;x,\xi)-x\xi-t\l, \ \varphi \in C^\infty((-T,T),S^{1,1}_\cl(\R^n))
$$
such that
\begin{itemize}
\item[$\bullet$]
$\partial_t\varphi(t;x,\xi)=q(x,d_x \varphi(t;x,\xi))$, $\varphi(0;x,\xi)=x\xi$;
\item[$\bullet$]
$C^{-1} \langle \xi \rangle\le \langle d_x \varphi(t;x,\xi)\rangle \le C
\langle \xi \rangle$, for a suitable constant $C>1$;
\item[$\bullet$]$q\in S_{{\rm{cl}}}^{m,1}(\R^n)$, $0<m<1$, $SG$-elliptic.
\end{itemize}
Since $q^{-1}(x,\xi)\in O(\norm{x}^{-1}\norm{\xi}^{-m})$
for $|x|+|\xi|\ge R>0$, it is not restrictive to assume that this estimate
holds on the whole phase space, so that,
for a certain constant $A>1$, 
\begin{equation}
	\label{eq:ellq}
	A^{-1}\norm{x}\norm{\xi}^m \le q(x,\xi) \le A\norm{x}\norm{\xi}^m.
\end{equation}
\begin{remark}
The assumption on $q^{-1}$ above amounts, at most, to modifying $q$ by adding and substracting
a compactly supported symbol, that is, an element of $S^{-\infty,-\infty}(\R^n)$. The
corresponding solutions $\varphi$ and $a$ of the eikonal and transport equations, respectively, would then change, at most, by an element of $C^\infty((-T,T),S^{-\infty,-\infty}(\R^n))$,
see \cite{Coriasco:998.2,Coriasco-Panarese:001,Coriasco:998.3}:
it is immediate, by integration by parts with respect to $t$, that an integral as \eqref{I} is 
$O(|\l|^{-\infty})$ for $a\in C^\infty((-T,T),S^{-\infty,-\infty}(\R^n))$. Then, the modified $q$ obviously keeps the same sign everywhere.
\end{remark}

For two functions $f,g$, defined on a common subset $X$ of $\R^{d_1}$ and depending on parameters
$y\in Y\subseteq\R^{d_2}$, we will write $f \prec g$ or $f(x,y) \prec g(x,y)$ to mean that there exists a suitable constant
$c>0$ such that $|f(x,y)|\le c |g(x,y)|$ for all $(x,y)\in X\times Y$. The notation $f\sim g$ or $f(x,y)\sim g(x,y)$
means that both $f\prec g$ and $g\prec f$ hold. 

\begin{remark}
\label{asymp-inf}
The ellipticity of $q$ yields, for $\l<0$, 
$$
\partial_t \phasef (t;x,\xi;\l)=q(x,d_x \varphi(t;x,\xi))-\l \succ
\langle x \rangle \langle \xi \rangle ^m +|\l|
$$ 
which, by integration by parts, implies 
$I(\l)={\rm{O}}(|\l|^{-\infty})$ when $\l \to - \infty$.
\end{remark}

\noindent
From now on any asimptotic estimate is to be meant for 
$\l \to +\infty$.

\vspace{3mm}

We will make use of a partition of unity on the phase space: the supports of
its elements will depend on suitably large positive constants $k_1,k_2>1$.
We also assume, as it is possible, $\l\ge\l_0$, again with an appropriate $\l_0>>1$.
As we will see below, the values of $k_1$, $k_2$ and $\l_0$ depend only on $q$ and its
associated seminorms. 

\begin{proposition}
\label{prop:h1}
Let $H_1$ be any function in $C^{\infty}_0(\R)$ such that 
${\rm{supp}}\,H_1 \subseteq[(2k_1)^{-1},2k_1]$, $0\leq H_1 \leq 1$ and 
$H_1\equiv1$ on $[k_1^{-1},k_1]$, where $k_1>1$ is a suitably chosen, large 
positive constant. Then  
\begin{equation}
\label{decomp1}
I(\l)={\rm{O}}(\l^{-\infty})+ 
\int e^{i\phasef(t;x,\xi;\l)}\psi(t)\,
H_1\!\left(\frac{\langle x \rangle \langle \xi \rangle^m}{\l}\right) 
a(t;x,\xi)\, dt  \dcut \xi  dx.
\end{equation}
\end{proposition}
\noindent
\begin{Pf}{.} Write 
\begin{equation}
\label{eq:split1}
\begin{aligned}
I(\l)&=\int e^{i\phasef(t;x,\xi;\l)}\psi(t)
\left[1- H_1\!\left(\frac{\langle x \rangle \langle \xi \rangle^m}{\l}\right)\right] 
a(t;x,\xi) \,dt \dcut \xi  dx
\\
&+\int e^{i\phasef(t;x,\xi;\l)}\psi(t)\,
H_1\!\left(\frac{\langle x \rangle \langle \xi \rangle^m}{\l}\right) 
a(t;x,\xi)\, dt  \dcut \xi  dx
\end{aligned}
\end{equation}
and observe that, by $ A^{-1}\langle x \rangle \langle \xi \rangle^m
\leq q(x,\xi)\leq A \langle x \rangle \langle \xi \rangle^m$, $x,\xi\in\R^n$, we find
\newlength{\esp}
\settowidth{\esp}{$^{-1}$}
\[
\begin{aligned}
|\partial_t \phasef(t;x,\xi;\l)|&\geq \frac{\l}{2}+\left(\frac{k_1}{2}-AC\right)
\langle x \rangle \langle \xi \rangle^m \qquad \qquad \qquad{\rm{when\ }}
\frac{\langle x \rangle \langle \xi \rangle^m}{\l}\leq k_1^{-1},
\\
|\partial_t \phasef(t;x,\xi;\l)|&\geq \frac{(AC)^{-1}}{2}\langle x \rangle \langle \xi \rangle^m
+\left[\frac{(AC)^{-1}}{2}k_1-1\right]\l \quad {\rm{when\ }}
\frac{\langle x \rangle \langle \xi \rangle^m}{\l}\geq k_1.\hspace{\esp}
\end{aligned}
\]
Thus, if $k_1 > 2AC$ we have
$|\partial_t \phasef(t;x,\xi;\l)| \sim \l+\langle x \rangle \langle \xi \rangle^m$
on the support of
$1- H_1\left(\dfrac{\langle x \rangle \langle \xi \rangle^m}{\l}\right)$, and the assertion follows
integrating by parts with respect to $t$ in the first integral of \eqref{eq:split1}.
\end{Pf}

\begin{remark}
	We actually choose $k_1>4AC>2AC$, since this will be needed in the proof of Proposition
	\ref{prop:i2} below, see also subsection \ref{sec:a1.4} in the Appendix.
\end{remark}

\noindent
Let us now pick $H_2\in C^\infty_0(\R)$ such that 
$0\leq H_2(\upsilon)\leq 1$, $H_2(\upsilon)=1$ for $|\upsilon|\leq k_2$ and $H_2(\upsilon)=0$ for 
$|\upsilon|\geq 2k_2$,  where $k_2>1$ is a constant which
we will choose big enough (see below). We can then write
\begin{align*}
(\l)={\rm{O}}(\l^{-\infty})&+ 
\int e^{i\phasef(t;x,\xi;\l)}\psi(t)
\, H_1\!\left(\frac{\langle x \rangle \langle \xi \rangle^m}{\l}\right) 
\, H_2(|\xi|) \, a(t;x,\xi)\, dt \dcut \xi  dx
\\
&+
\int e^{i\phasef(t;x,\xi;\l)}\psi(t)
\, H_1\!\left(\frac{\langle x \rangle \langle \xi \rangle^m}{\l}\right) 
[1-H_2(|\xi|)] \, a(t;x,\xi)\, dt \dcut \xi dx
\\
=  {\rm{O}}(\l^{-\infty})&+ I_1(\l) + I_2(\l).
\nonumber
\end{align*}

In what follows, we will sistematically use the notation $S^{r,\rho}=S^{r,\rho}(y,\eta)$, $y\in\R^k$, $\eta\in\R^l$,
to generally denote functions depending smoothly on $y$ and $\eta$ and satisfying $SG$-type estimates of order $r,\rho$ in $y,\eta$.
In a similar fashion, $S_T^{r,\rho}=C^\infty((-T,T),S^{r,\rho}(y,\eta))$ will stand for some function of the same kind which, additionally, depends smoothly on $t\in (-T,T)$, 
and, for all $s\in\Z_+$,
$D^s_t C^\infty((-T,T),S^{r,\rho}(y,\eta))$ satisfies $SG$-type estimates of order $r,\rho$ in $y,\eta$, uniformly with respect to $t\in (-T,T)$. 

To estimate $I_1(\l)$, we will apply
 the Stationary Phase Theorem.
We begin by rewriting the integral $I_1(\l)$, using the fact that $\varphi$ is solution of the eikonal equation associated with $q$ and that $q$ is a classical $SG$-symbol. 
Note that then $\partial^2_t \varphi \in C^\infty((-T,T), S_{{\rm{cl}}}^{2m-1,1}(\R^n))\subseteq C^\infty((-T,T),S_{{\rm{cl}}}^{m,1}(\R^n))$, since
$$
\partial^2_t \varphi (t;x,\xi)=\sum_{i=1}^{n}
(\partial_{\xi{_i}}q)(x,d_x \varphi(t;x,\xi)) \, \partial_{x_{i}}\big(q(x,d_x \varphi(t;x,\xi))\big).
$$

In view of
the Taylor expansion of $\varphi$ at $t=0$, recalling the property 
$q(x,\xi)=\omega(x) q_e(x,\xi)+S^{m,0}(x,\xi)$, $\omega$ a fixed $0$-excision function,
we have, for some $0<\delta_1<1$,
\begin{align*}
	\phasef(t;x,\xi;\l)&=-\l t -x\xi + \varphi(0;x,\xi)+t\,\partial_t\varphi(0;x,\xi)+
	\frac{t^2}{2}\,\partial_t^2\varphi(t \delta_1;x,\xi)
	\\
	&=-\l t + tq(x,\xi)+t^2S_T^{2m-1,1}(x,\xi)
	\\
	&=-\l t + t\omega(x) q_e(x,\xi) +tS^{m,0}(x,\xi)+ 
	t^2 S_T^{2m-1,1}(x,\xi)
	\\
	&=-\l t + t\omega(x) q_e(x,\xi) +tS^{m,0}(x,\xi)+ 
	t^2 \omega(x)S_{T,e}^{2m-1,1}(x,\xi)+t^2S_T^{2m-1,0}(x,\xi),
\end{align*}
where the subscript $e$ denotes the $x$-homogeneous (exit) principal parts of the
involved symbols, which are all $SG$-classical and real-valued, see \cite{Coriasco-Panarese:001}.

Observe that $|x|\sim \l$ on the support of the integrand in $I_1(\l)$, so that we can, in fact, assume $\omega(x)\equiv 1$ there. Indeed, recalling that, by definition, $\omega\in C^\infty(\R^n)$, $\omega(\upsilon)\equiv 0$ for $|\upsilon|\le B$, $\omega(\upsilon)\equiv 1$ for $|\upsilon|\ge 2B$, with a fixed constant $B>0$, it is enough to observe that
\[
|\xi|\prec 1,\norm{x}\norm{\xi}^m\sim\l\Rightarrow\norm{x}\sim\l,
\]
which of course implies $\norm{x}\sim|x|$, provided $\l_0\le\l$ is
large enough. Moreover, by the ellipticity of $q$, writing $x=|x|\varsigma$, $\varsigma\in\Sn$,
with the constant $A>1$ of \eqref{eq:ellq},
\begin{align}
\nonumber
A^{-1}\norm{x}\norm{\xi}^m&\le q(x,\xi)=\omega(x)q_e(x,\xi)+S^{m,0}(x,\xi)\le A \norm{x}\norm{\xi}^m
\\
\nonumber
&\Rightarrow
A^{-1}\frac{\norm{x}}{|x|}\norm{\xi}^m
\le 
\omega(x) q_e(\varsigma,\xi)+\frac{S^{m,0}(x,\xi)}{|x|}
\le 
A \frac{\norm{x}}{|x|}\norm{\xi}^m
\\
\label{eq:ellqe}
&\Rightarrow
A^{-1}\norm{\xi}^m
\le 
q_e(\varsigma,\xi)
\le 
A \norm{\xi}^m,\quad\varsigma\in\Sn,\xi\in\R^n
\end{align}
taking the limit for $|x|\to+\infty$. 
Then, setting $x=\l\zeta\varsigma$, $\zeta\in[0,+\infty)$, $\varsigma\in \Sn$,
$\l\ge\l_0>>1$, in $I_1(\l)$, by homogeneity and the above remarks, we can write
\begin{align*}
	\phasef(t;\l\zeta\varsigma,\xi;\l)&=-\l t + t\omega(\l\zeta\varsigma) q_e(\l\zeta\varsigma,\xi) +tS^{m,0}(\l\zeta\varsigma,\xi)
	\\
	&+ 
	t^2 \omega(\l\zeta\varsigma)S_{T,e}^{2m-1,1}(\l\zeta\varsigma,\xi)+t^2S_T^{2m-1,0}(\l\zeta\varsigma,\xi)
	\\
	&=-\l t + \l\zeta t q_e(\varsigma,\xi) + 
	\l\zeta t^2  S_{T,e}^{2m-1,1}(\varsigma,\xi)+
	tS^{m,0}(\l\zeta\varsigma,\xi)+t^2S_T^{2m-1,0}(\l\zeta\varsigma,\xi)
	\\
	&=\l [-t + \zeta t q_e(\varsigma,\xi) + \zeta t^2  S_{T,e}^{2m-1,1}(\varsigma,\xi)]+
	G_1(\l;t,\zeta;\varsigma,\xi)
	\\
	&=\l F_1(t,\zeta;\varsigma,\xi)+G_1(\l;t,\zeta;\varsigma,\xi),
\end{align*}
and find, in view of the compactness of the support of the integrand (see the proof of Proposition \ref{prop:i1}
below) and the hypotheses,
\begin{align*}
I_1(\l)  &= \l^n \int e^{i\l F_1(t,\zeta;\varsigma,\xi)}
\,e^{i G_1(\l;t,\zeta;\varsigma,\xi)} \, \psi(t) \, 
a(t;\l\zeta\varsigma,\xi)\,  H_1\!\left(\frac{\norm{\l\zeta}\norm{\xi}^m}{\l}\right) H_2(|\xi|) \zeta^{n-1}\, dt d\zeta \dcut \xi d\varsigma
\\& = \frac{\l^n}{(2\pi)^n} \int e^{i\l F_1(X,Y)} U_1(X,Y;\l)\,dXdY,
\end{align*}
with $X=(t,\zeta)$, $Y=(\varsigma,\xi)$. We can now prove
\begin{proposition} 
\label{prop:i1}
Choosing the constants $k_1,\l_0>1$ large enough and $T>0$ suitably small,
we have, for any $k_2>1$ and for a certain sequence $c_j$, $j=0,1,\dots$,
\[
	I_1(\l)\sim \sum_{j=0}^{+\infty}c_j \l^{n-1-j},
\]
that is, $I_1(\l)=c_0 \l^{n-1}+O(\l^{n-2})$, with
\[
c_0=\frac{1}{(2\pi)^{n-1}}\int_{\R^n}\int_{\Sn}\frac{H_2(|\xi|)}{q_e(\varsigma,\xi)^{n}}d\varsigma d\xi.
\]
\end{proposition}
\noindent
\begin{Pf}{.}
It is easy to see that, on the support
of $U_1$, the phase function
$F_1(X,Y)$ admits a unique, nondegenerate, stationary point 
$X_0=X_0(Y)=(0,q_e(\varsigma,\xi)^{-1})$, i.e. $F^\prime_{1,X}(X_0(Y),Y)=0$ for all $Y$
such that $(X,Y)\in\supp\,U_1$, provided $T>0$ is chosen
suitably small (see, e.g., \cite{GriSjo:994}, p. 136), and the Hessian $\det(F_{1,X}^{\prime\prime}(X_0(Y),Y))$ equals $-q_e(\varsigma,\xi)^2<0$. Moreover, the amplitude function
\[
	U_1(X,Y;\l)=
	\psi(t) \, 
  H_1\!\left(\frac{\norm{\l\zeta}\norm{\xi}^m}{\l}\right) H_2(|\xi|) \,
  a(t;\l\zeta\varsigma,\xi)\,\zeta^{n-1}\,e^{i G(\l;t,\zeta;\varsigma,\xi)}
\]
is compactly supported with respect to the variables $X$ and $Y$, and satisfies, for all
$\gamma\in\Z_+^2$,
\[
D^\gamma_X U_1(X,Y;\l) \prec 1
\]
for all $X$, $Y$, $\l\ge\l_0$. In fact:
\begin{enumerate}
	\item $\psi\in C_0^\infty((-T,T))$, $\varsigma\in\Sn$, 
	$\supp [H_2(|\xi|)]\subseteq\{\xi\colon|\xi|\le 2k_2\}$, and
	\[
		(2k_1)^{-1}\le \norm{\xi}^m\sqrt{\frac{1}{\l^2}+\zeta^2}\le 2k_1
		\Rightarrow 0<\sqrt{\frac{1}{4k_1^2\norm{2k_2}^{2m}}-\frac{1}{\l_0^2}}\le 
		\zeta\le2k_1,
	\]
	where $\l_0>2k_1\norm{2k_2}^m$;
	\item all the factors appearing in the expression of $U_1$ are uniformly bounded, 
	together with all their $X$-derivatives, for $X\in S_X=\supp\,\psi\times[\zeta_0,\zeta_1]$,
	$Y\in S_Y=\Sn\times\{\xi\colon|\xi|\le 2k_2\}$, and $\l\ge\l_0$.
\end{enumerate}
Of course, (2) trivially holds
for the cutoff functions $\psi(t)$ and $H_2(|\xi|)$, and for the factor $\zeta^{n-1}$.
Since $a(t;x,\xi)\in S_T^{0,0}(x,\xi)$, on $S_X\times S_Y$ we have, for all $\gamma\in\Z_+^2$ and $\l\ge\l_0>1$,
\[
	D^\gamma_X a(t;\l\zeta\varsigma,\xi)\prec 
	\norm{\l\zeta}^{-\gamma_2}\l^{\gamma_2}\norm{\xi}^m
	\prec \frac{1}{\left(\dfrac{1}{\l^2}+\zeta^2\right)^{\frac{\gamma_2}{2}}}
	<\frac{1}{\zeta^{\gamma_2}}\prec 1.
\]
Moreover, since $G_1\in S_T^{m,0}(x,\xi)$ is actually in $S_T^{-\infty,0}(x,\xi)\subset S_T^{0,0}(x,\xi)$ on $S_X\times S_Y$,  the same holds for $\exp(iG_1)$, by an application of the Fa\`{a} di Bruno formula for the derivatives of compositions of functions, so also this factor fulfills the desired estimates. Finally, another straightforward computation shows that, for all $\gamma_2\in\Z_+$ and $\l\ge\l_0>1$, 
\[
	D_\zeta^{\gamma_2} H_1\!\left(\frac{\norm{\l\zeta}\norm{\xi}^m}{\l}\right) \prec 1
\]
on $S_X\times S_Y$.
The Proposition is then a consequence of the Stationary Phase Theorem (see \cite{D96},
Proposition 1.2.4, \cite{Ho1}, Theorem 7.7.6), applied to the 
integral with respect to $X=(t,\zeta)$. In particular, the leading term is given by
$\dfrac{\l^n}{(2\pi)^{n-1}}$
times the integral with respect to $Y$ of $\l^{-1}|\det(F_{1,X}^{\prime\prime}(X_0(Y),Y))|^{-\frac{1}{2}} \, U_1(X_0(Y),Y;\l)$, that is
\begin{align*}
I_1(\l)&=\frac{\l^{n-1}}{(2\pi)^{n-1}}\int_{\R^n}\int_{\Sn}
\frac{1}{q_e(\varsigma,\xi)} \,\psi(0)\,H_1\!\left(\frac{\norm{\frac{\l}{q_e(\varsigma,\xi)}}\norm{\xi}^m}{\l}\right)
\frac{H_2(|\xi|)}{q_e(\varsigma,\xi)^{n-1}}
\,a\!\left(0;\frac{\l\varsigma}{q_e(\varsigma,\xi)},\xi\right)
d\varsigma d\xi
\\
&+O(\l^{n-2})
\\
&=\frac{\l^{n-1}}{(2\pi)^{n-1}}\int_{\R^n}\int_{\Sn}
H_1\!\left(\frac{\norm{\frac{\l}{q_e(\varsigma,\xi)}}\norm{\xi}^m}{\l}\right)
\frac{H_2(|\xi|)}{q_e(\varsigma,\xi)^{n}}
d\varsigma d\xi+O(\l^{n-2})
\\
&=\frac{\l^{n-1}}{(2\pi)^{n-1}}\int_{\R^n}\int_{\Sn}
\frac{H_2(|\xi|)}{q_e(\varsigma,\xi)^{n}}
d\varsigma d\xi+O(\l^{n-2}).
\end{align*}
recalling that $\psi(0)=1$, $a(0;x,\xi)=1$ for all $x,\xi\in\R^n$. Indeed, having chosen $k_1>2A$, $\l_0>2k_1\norm{2k_2}^m$, \eqref{eq:ellqe} implies
\[
k_1^{-1}<A^{-1}<
\frac{\norm{\frac{\l}{q_e(\varsigma,\xi)}}\norm{\xi}^m}{\l}=
\sqrt{\left(\frac{\norm{\xi}^{m}}{\l}\right)^2+\left(\frac{\norm{\xi}^m}{q_e(\varsigma,\xi)}\right)^2}
<\sqrt{\frac{1}{4k_1^2}+A^2}< k_1,
\]
uniformly for $\varsigma\in\Sn$, $\xi\in\supp[H_2(|\xi|)]$, $\l\ge\l_0$. This concludes the proof. 
\end{Pf}

\vspace{3mm}

Let us now consider $I_2(\l)$. We follow a procedure close to that used in the 
proof of Theorem 7.7.6 of \cite{Ho1}. However, since here we lack
the compactness of the support of the amplitude with respect to $x$, we need
explicit estimates to show that all the involved integrals are
convergent, so we give below the argument in full detail.

We initially proceed as in the analysis of $I_1(\l)$ above.
In view of the presence of the factor $1-H_2(|\xi|)$
in the integrand, we can now assume $|\xi|\ge k_2>\max\{B,1\}$, 
$B>0$ the radius of the smallest ball in $\R^n$ including 
$\supp\,(1-\omega)$, so that $q(x,\xi)=\omega(\xi)q_\psi(x,\xi)+S^{m-1,1}(x,\xi)=q_\psi(x,\xi)+S^{m-1,1}(x,\xi)$. Then, with some $0<\delta_2<1$,
\begin{align*}
	\phasef(t;x,\xi;\l)&=-\l t -x\xi + \varphi(0;x,\xi)+t\,\partial_t\varphi(0;x,\xi)+
	\frac{t^2}{2}\,\partial_t^2\varphi(t \delta_2;x,\xi)
	\\
	&=-\l t + tq(x,\xi)+t^2S_T^{2m-1,1}(x,\xi)
	\\
	&=-\l t + t q_\psi(x,\xi) +tS^{m-1,1}(x,\xi)+ 
	t^2 S_T^{2m-1,1}(x,\xi).
\end{align*}
Setting $\xi=(\l\zeta)^\frac{1}{m}\varsigma$, $\zeta\in[0,+\infty)$, $\varsigma\in \Sn$,
$\l\ge\l_0$, we can rewrite $I_2(\l)$ as
\begin{align*}
I_2(\l)  &= \frac{n}{m}\frac{\l^\frac{n}{m}}{(2\pi)^n} \int e^{i \l(- t + \zeta t q_\psi(x,\varsigma) +t\l^{-1}S^{m-1,1}(x,(\l\zeta)^\frac{1}{m}\varsigma)+ 
	t^2 \l^{-1}S_T^{2m-1,1}(x,(\l\zeta)^\frac{1}{m}\varsigma))} \cdot
	\\ 
	&\phantom{\frac{n}{m}\frac{\l^\frac{n}{m}}{(2\pi)^n} \int}\cdot\psi(t) \, 
a(t;x,(\l\zeta)^\frac{1}{m}\varsigma)\,  H_1\!\left(\frac{\norm{x}\norm{(\l\zeta)^\frac{1}{m}\varsigma}^m}{\l}\right) \left[1-H_2((\l\zeta)^\frac{1}{m})\right] \, \zeta^{\frac{n}{m}-1}\, dt d\zeta d\varsigma dx
\\& = \frac{n}{m}\frac{\l^\frac{n}{m}}{(2\pi)^n}\int e^{i\l F_2(X,Y;\l)} U_2(X,Y;\l)\,dXdY,
\end{align*}
$X=(t,\zeta)$, $Y=(\varsigma,x)$, where we have set
\begin{align*}
	F_2(X,Y;\l)&= - t + \zeta t q_\psi(x,\varsigma) +
	t\l^{-1}S^{m-1,1}(x,(\l\zeta)^\frac{1}{m}\varsigma)+ 
	t^2 \l^{-1}S_T^{2m-1,1}(x,(\l\zeta)^\frac{1}{m}\varsigma)
	\\
	U_2(X,Y;\l)&=\psi(t) \,   
	H_1\!\left(\frac{\norm{x}\norm{(\l\zeta)^\frac{1}{m}\varsigma}^m}{\l}\right) 
	\left[1-H_2((\l\zeta)^\frac{1}{m})\right] \, a(t;x,(\l\zeta)^\frac{1}{m}\varsigma)\,\zeta^{\frac{n}{m}-1}.
\end{align*}
On the support of $U_2$ we have
\[
	\frac{\norm{x}\norm{(\l\zeta)^\frac{1}{m}\varsigma}^m}{\l}\sim 1
	\;\mbox{  and  }\;
	(\l\zeta)^\frac{1}{m}\succ 1\Rightarrow 
	\norm{(\l\zeta)^\frac{1}{m}\varsigma}^m = \norm{(\l\zeta)^\frac{1}{m}}^m
	\sim \l\zeta,
\]
so that
\begin{equation}
	\label{eq:suppu2}
	\frac{\norm{x}\l\zeta}{\l}\sim 1\Leftrightarrow \zeta \sim \norm{x}^{-1}
	\mbox{ and } |x|<\norm{x}\le 2k_1\norm{k_2}^{-m} \l = \tilde{\varkappa}\l.
\end{equation}
For any fixed $Y\in\Sn\times\R^n$ we then have $X$ belonging to
a compact set, uniformly with respect to $\l\ge\l_0$, say 
$\supp\, \psi \times [c^{-1}\norm{x}^{-1},c\norm{x}^{-1}]$, for a suitable $c>1$.

\begin{remark}
\label{rem:esti2}
Incidentally, we observe that a rough estimate of $\l^\frac{n}{m}I_2(\l)$ is
\begin{align*}
	\int e^{i\l F_2(X,Y;\l)} U_2(X,Y;\l)\,dX \prec 
	\,&\norm{x}^{-\frac{n}{m}+1}\int_{c^{-1}\norm{x}^{-1}}^{c\norm{x}^{-1}} d\zeta\prec
	\norm{x}^{-\frac{n}{m}}
	\\
	&\Rightarrow \l^\frac{n}{m}\int e^{i\l F_2(X,Y;\l)} U_2(X,Y;\l)\,dXdY\prec\l^n,\l\to+\infty.
\end{align*}
An even less precise result would be the bound $\l^\frac{n}{m}$, using the convergence
of the integral with respect to $x$ in the whole $\R^n$, given by $-\dfrac{n}{m}+n<0$.
\end{remark}

\noindent
The next Lemma is immediate, and we omit the proof:

\begin{lemma}
	\label{lemma:simbcomp}
	$S^{s,\sigma}_T(x,(\l\zeta)^\frac{1}{m}\varsigma)
	=S^{s,\sigma}_T(x,(\l\zeta)^\frac{1}{m})$
	for any $\zeta\in[0,+\infty)$, $x\in\R^n$, $\varsigma\in\Sn$,
	$\l\ge\l_0$, $m\in(0,1)$, and, for all $\gamma\in\Z_+^2$,
	$$
	D^\gamma_X S^{s,\sigma}_T(x,(\l\zeta)^\frac{1}{m}) = 
	\zeta^{-\gamma_2} S^{s,\sigma}_T(x,(\l\zeta)^\frac{1}{m}).
	$$
\end{lemma}

\noindent
The main result of this Section is

\begin{proposition}
\label{prop:i2}
If $k_1,k_2,\l_0>1$ are chosen large enough we have
\begin{equation}
	\label{eq:i2weyl}
	I_2(\l)=\frac{n}{m}\,d_0 \l^{\frac{n}{m}-1}+O(\l^{n-1})+O(\l^{\frac{n}{m}-2}).
\end{equation}
Explicitely,
\[
d_0=
\frac{1}{(2\pi)^{n-1}}\int_{\R^n}\int_{\Sn}
\frac{1}{q_\psi(x,\varsigma)^\frac{n}{m}}\,d\varsigma dx.
\]
\end{proposition}

\vspace{5mm}

We will prove Proposition \ref{prop:i2} through various intermediate steps. First of all,
arguing as in the proof of \eqref{eq:ellqe}, exchanging the role of $x$ and $\xi$,
we note that, for all $x\in\R^n$, $\varsigma\in\Sn$,
\begin{equation}
	\label{eq:ellqpsi}
	A^{-1}\norm{x}\le q_\psi(x,\varsigma)\le A\norm{x},
\end{equation}
$(x,\varsigma)\in\R^n\times\Sn$. We now study 
\begin{align*}
F_{2,X}^\prime(X,Y;\l)&=\left(
\begin{array}{c}
\partial_t F_{2}(X,Y;\l)
\\
\partial_\zeta F_{2}(X,Y;\l)
\end{array}
\right)
\\
&=
\left(
		\begin{array}{l}
			-1+\dfrac{\zeta}{\zeta_0}+
			\l^{-1}S^{m-1,1}(x,(\l\zeta)^\frac{1}{m})+ 
			t \l^{-1}S_T^{2m-1,1}(x,(\l\zeta)^\frac{1}{m})
		\\
			t(q_\psi(x,\varsigma) +
			\l^{-1}\zeta^{-1}S^{m-1,1}(x,(\l\zeta)^\frac{1}{m})+ 
			t \l^{-1}\zeta^{-1}S_T^{2m-1,1}(x,(\l\zeta)^\frac{1}{m}))
		\end{array}
\right),
\end{align*}
$X=(t,\zeta)\in S_X=\supp\,\psi\times [c^{-1}\norm{x}^{-1},c\norm{x}^{-1}]$,
$Y=(\varsigma,x)\in S_Y=\Sn\times\R^n$, $\l\ge\l_0$, 
where we have used Lemma \ref{lemma:simbcomp}.
By the symbolic calculus, remembering that $\l\zeta\ge k_2^m>1$ on $\supp\,U_2$,
we can rewrite the expressions above as
\begin{align*}
		\partial_t F_{2}(X,Y;\l)&=-1+\dfrac{\zeta}{\zeta_0}+
			\zeta(\l\zeta)^{-1}S^{m-1,1}(x,(\l\zeta)^\frac{1}{m})+ 
			t \zeta(\l\zeta)^{-1}S_T^{2m-1,1}(x,(\l\zeta)^\frac{1}{m})
			\\
			&=-1+\dfrac{\zeta}{\zeta_0}+
			\zeta[(\l\zeta)^\frac{1}{m}]^{-m}S^{m-1,1}(x,(\l\zeta)^\frac{1}{m})+ 
			t \zeta[(\l\zeta)^\frac{1}{m}]^{-m}S_T^{2m-1,1}(x,(\l\zeta)^\frac{1}{m})
			\\
			&=-1+\dfrac{\zeta}{\zeta_0}+
			\zeta S^{-1,1}(x,(\l\zeta)^\frac{1}{m})+ 
			t \zeta S_T^{m-1,1}(x,(\l\zeta)^\frac{1}{m}),
		\\
		\partial_\zeta F_{2}(X,Y;\l)&=t(q_\psi(x,\varsigma) +
			S^{-1,1}(x,(\l\zeta)^\frac{1}{m})+ 
			t S_T^{m-1,1}(x,(\l\zeta)^\frac{1}{m})).
\end{align*}
It is clear that $\zeta\sim\norm{x}^{-1}$ implies
$\zeta S^{-1,1}(x,(\l\zeta)^\frac{1}{m})=S^{-1,0}(x,(\l\zeta)^\frac{1}{m})$ and
$\zeta S_T^{m-1,1}(x,(\l\zeta)^\frac{1}{m})=S_T^{m-1,0}(x,(\l\zeta)^\frac{1}{m})$,
so that we finally have
\begin{align*}
			\partial_t F_2(X,Y;\l)&=-1+\dfrac{\zeta}{\zeta_0}+
			S^{-1,0}(x,(\l\zeta)^\frac{1}{m})+ 
			t S_T^{m-1,0}(x,(\l\zeta)^\frac{1}{m}),
		\\
			\partial_\zeta F_2(X,Y;\l)&=t(q_\psi(x,\varsigma) +
			S^{-1,1}(x,(\l\zeta)^\frac{1}{m})+ 
			t S_T^{m-1,1}(x,(\l\zeta)^\frac{1}{m})).
\end{align*}
We now prove that, modulo an $O(|\l|^{-\infty})$ term, we can consider an amplitude such
that, on its support, the ration $\zeta/\zeta_0$ is very close to $1$. To this aim,
take $H_3\in C^\infty_0(\R)$ such that 
$0\leq H_3(\upsilon)\leq 1$, $H_3(\upsilon)=1$ for $|\upsilon|\leq \dfrac{3}{2}\varepsilon$
and $H_3(\upsilon)=0$ for $|\upsilon|\geq 2\varepsilon$, with
an arbitrarily fixed, small enough $\varepsilon\in\left(0,\dfrac{1}{2}\right)$, and set
\[
	V_1(X,Y;\l)=U_2(X,Y;\l)\cdot\left[1-H_3\left(\dfrac{\zeta}{\zeta_0}-1\right)\right],
	\;
	V_2(X,Y;\l)=U_2(X,Y;\l)\cdot H_3\left(\dfrac{\zeta}{\zeta_0}-1\right),
\]
\[
	J_1(\l)=\int e^{i\l F_2(X,Y;\l)}V_1(X,Y;\l)\,dXdY,
	\;
	J_2(\l)=\int e^{i\l F_2(X,Y;\l)}V_2(X,Y;\l)\,dXdY.
\]
\begin{proposition}
	\label{prop:j1}
	With the choices of $T, k_1,\l_0$ above, for any $\varepsilon\in\left(0,\dfrac{1}{2}\right)$
	we can find $k_2>1$ large enough such that $J_1(\l)=O(\l^{-\infty})$.
\end{proposition}
\noindent
\begin{Pf}{.}
	Since $0<m<1$,
	in view of \eqref{disSG}, \eqref{eq:suppu2}, and \eqref{eq:ellqpsi}, 
	we can choose $k_2>1$ so large that,
	for an arbitrarily fixed $\varepsilon\in\left(0,\dfrac{1}{2}\right)$, for any 
	$\l\ge\l_0$, $\zeta\in(0,+\infty)$ satisfying $|\xi|=(\l\zeta)^\frac{1}{m}\ge k_2$,
	\begin{equation}
		\label{eq:gradf2}
		\begin{aligned}
		\mbox{in $\partial_t F_2(X,Y;\l)$, }&
		\left|S^{-1,0}(x,(\l\zeta)^\frac{1}{m})\right|\le\dfrac{\varepsilon}{2},
		\quad
		\left|t S_T^{m-1,0}(x,(\l\zeta)^\frac{1}{m})\right|\le\dfrac{\varepsilon}{2},
		\\
		\mbox{ and }&
		\left|\zeta_0\frac{d}{d\zeta}S^{-1,0}(x,(\l\zeta)^\frac{1}{m})\right|
		=\left|\zeta_0\zeta^{-1} S^{-1,0}(x,(\l\zeta)^\frac{1}{m})\right|\le k_0<1,
		\\
		\mbox{in $\partial_\zeta F_2(X,Y;\l)$, }&
		|S^{-1,1}(x,(\l\zeta)^\frac{1}{m})+ 
			t S_T^{m-1,1}(x,(\l\zeta)^\frac{1}{m}))|\le\frac{A^{-1}}{2}\norm{x},
		\end{aligned}
	\end{equation}
	uniformly with respect to $(X,Y)\in S_X\times S_Y\supseteq\supp\,U_2(.;\l)$.
	Then, $F_2$ is non-stationary on $\supp\,V_1$, since
	there we have 
	$\left|\dfrac{\zeta}{\zeta_0}-1\right|\ge\dfrac{3}{2}\varepsilon$, while 
	\[
	\left|S^{-1,0}(x,(\l\zeta)^\frac{1}{m})+
	t S_T^{m-1,0}(x,(\l\zeta)^\frac{1}{m})\right|\le\varepsilon,
	\]
	which implies
	$\partial_t F_2(X,Y;\l)\succ 1$.
	Observing that, on $\supp\,V_1$,
	$\partial_t F_2(X,Y;\l) = S_T^{0,0}(x,(\l\zeta)^\frac{1}{m})$, as well as $V_1(X,Y;\l)
	= S_T^{0,0}(x,(\l\zeta)^\frac{1}{m})$,
	the assertion follows by repeated integrations by parts with respect to $t$,
	using the operator
	\[
		L_t=\dfrac{1}{\l\partial_t F_2(X,Y;\l)}D_t \Rightarrow L_1 e^{i\l F_2(X,Y;\l)}=e^{i\l F_2(X,Y;\l)}
	\]
	and recalling Remark \ref{rem:esti2}.
\end{Pf}

\vspace{3mm}

\begin{proposition}
With the choices of $\varepsilon, T>0$, $k_1,k_2,\l_0>1$ above,
we can assume, modulo an $O(\l^{n-1})$ term, that the integral with respect to
$x$ in $J_2(\l)$ is extended to the set $\{x\in\R^n\colon\norm{x}\le \varkappa\l\}$, with 
\begin{equation}
	\label{eq:varkappa}
	\varkappa=\left(1-\frac{\varepsilon}{2}\right)[A(2k_2)^m]^{-1}.
\end{equation}
\end{proposition}
\noindent
\begin{Pf}{.}
Indeed, if $\varkappa < \tilde{\varkappa}=2k_1\norm{k_2}^{-m}$, we can
split $J_2(\l)$ into the sum
\begin{equation}
	\label{eq:j2split}
	\int_{\varkappa\l\le\norm{x}\le\tilde{\varkappa}\l}\int_{\Sn}\int e^{i\l F_2}V_2\,dX d\varsigma dx
	+
	\int_{\norm{x}\le\varkappa\l}\int_{\Sn}\int e^{i\l F_2}V_2\,dX d\varsigma dx,
\end{equation}
since the inequality $\varkappa < \tilde{\varkappa}$ is true when
$k_2$ is sufficiently large. Observing that, on $\supp\, U_2$,
\[
	\norm{x}\sim\l\Rightarrow \norm{\xi}^{m}=\frac{\norm{x}\norm{\xi}^m}{\l}\,\frac{\l}{\norm{x}}\sim1
	\Rightarrow |\xi|\le k_3,
\]
switching back to the original variables, 
the first integral in \eqref{eq:j2split} can be treated as $I_1(\l)$, and gives,
in view of Proposition \ref{prop:i1}, an $O(\l^{n-1})$ term, as stated.
\end{Pf}

\vspace{3mm}

Now we can show that $F_2(X,Y;\l)$ admits a unique, nondegenerate 
stationary point $X_0^*=X_0^*(Y,\l)$ belonging to $\supp\,V_2$ for $\norm{x}\le\varkappa\l$.
Under the same hypotheses, $X_0^*$ lies in a circular neighbourhood of
$X_0=(0,\zeta_0)=(0,q_\psi(x,\varsigma)^{-1})$ of arbitrarily small radius:
\begin{proposition}
	\label{prop:statF2}
	With $\varepsilon\in\left(0,\dfrac{1}{2}\right)$, $T>0$, $k_1,k_2,\l_0>1$ fixed above,
	$F_{2,X}^\prime(X,Y;\l)$ vanishes on
	$supp\,V_2$ only for $X=X_0^*(Y;\l)=(0,\zeta_0^*(Y;\l))$, 
	i.e., $F_{2,X}^\prime(X_0^*(Y;\l),Y;\l)=0$
	for all $Y$ such that $(X,Y;\l)\in\supp\,V_2$. Moreover,	
	\[
		\det(F^{\prime\prime}_{2,X}(X_0^*(Y;\l),Y))\sim\norm{x}^2 \mbox{ and }
		|X_0^*(Y;\l)-X_0(Y)|=|\zeta_0^*(Y;\l)-\zeta_0(Y)|\le \frac{A\varepsilon}{2}\norm{x}^{-1}
	\]
	holds on $\supp\,V_2$.
\end{proposition}
\noindent
\begin{Pf}{.} We have to solve
\[
\begin{cases}
			0&\hspace{-3mm}=-1+\dfrac{\zeta}{\zeta_0}+
			S^{-1,0}(x,(\l\zeta)^\frac{1}{m})+ 
			t S_T^{m-1,0}(x,(\l\zeta)^\frac{1}{m})
		\\
			0&\hspace{-3mm}=t(q_\psi(x,\varsigma) +
			S^{-1,1}(x,(\l\zeta)^\frac{1}{m})+ 
			t S_T^{m-1,1}(x,(\l\zeta)^\frac{1}{m})),
\end{cases}
\]
$(X,Y;\l)\in\supp\,V_2$. By \eqref{eq:ellqpsi} and \eqref{eq:gradf2}, with the choices
of $\varepsilon, T>0$, $k_1,k_2,\l_0$ above, the coefficient of $t$ in the second equation
does not vanish at any point of $\supp\,V_2$
 Then $t=0$,
	and $\zeta$ must satisfy
	\begin{equation}
		\label{eq:zeta0star}
		-1+\frac{\zeta}{\zeta_0}+S^{-1,0}(x,(\lambda\zeta)^\frac{1}{m})=0
		\Leftrightarrow
		\zeta=\zeta_0(1+S^{-1,0}(x,(\lambda\zeta)^\frac{1}{m})=G(\zeta;Y;\l).
	\end{equation}
	Since, by the choice of $k_2$, $|\partial_\zeta G(\zeta;Y;\l)|\le k_0<1$,
	uniformly with respect to $Y\in\Sn\times\{x\in\R^n\colon\norm{x}\le\varkappa\l\}$,
	$\l\ge\l_0$, $G$ has a unique
	fixed point $\zeta_0^*=\zeta_0^*(Y;\l)$, smoothly depending on the parameters,
	see the Appendix for more details. Since
	\begin{equation}
		\label{eq:f2second}
		\begin{aligned}
		\partial^2_t F_2(X,Y;\l)&=S_T^{m-1,0}(x,(\l\zeta)^\frac{1}{m}),
		\\
		\partial_t\partial_\zeta F_2(X,Y;\l)&=q_\psi(x,\varsigma)(1+\zeta_0\zeta^{-1}(
				S^{-1,0}(x,(\l\zeta)^\frac{1}{m})+ 
				t S_T^{m-1,0}(x,(\l\zeta)^\frac{1}{m}))),
		\\
		\partial_\zeta^2 F_2(X,Y;\l)&=t\zeta^{-1}(
				S^{-1,1}(x,(\l\zeta)^\frac{1}{m})+ 
				t S_T^{m-1,1}(x,(\l\zeta)^\frac{1}{m})),
		\end{aligned}
	\end{equation}
	we can assume that $\l\zeta\ge k_2^m$ and the choices of the
	other parameters imply, on $\supp\,V_2$,
	\[
		\partial^2_t F_2(X,Y;\l)\prec\frac{\varepsilon}{2},
		\quad
		\partial_t\partial_\zeta F_2(X,Y;\l)\sim\norm{x},
		\quad
		\partial^2_\zeta F_2(X,Y;\l)\prec\frac{\varepsilon}{2}\norm{x}^2.
	\]
	So we have proved that, on $\supp\,V_2$,
	\begin{equation}
		\label{eq:M}
		\begin{aligned}
		M&=F^{\prime\prime}_{2,X}(X_0^*(Y;\l),Y;\l)=
		\left(
			\begin{array}{cc}
				M_{11}
				&
				M_{12}
				\\
				M_{12}
				&
				0
			\end{array}
		\right)=
		\\&=
		\left(
			\begin{array}{cc}
				S_T^{m-1,0}(x,(\l\zeta_0^*)^\frac{1}{m})
				&
				q_\psi(x,\varsigma)\left[1+\dfrac{\zeta_0}{\zeta_0^*}
				S^{-1,0}(x,(\l\zeta_0^*)^\frac{1}{m})\right]
				\\
				q_\psi(x,\varsigma)\left[1+\dfrac{\zeta_0}{\zeta_0^*}
				S^{-1,0}(x,(\l\zeta_0^*)^\frac{1}{m})\right]
				&
				0
			\end{array}
		\right)
		\\
		&\Rightarrow \det(M)=
		-q_\psi(x,\varsigma)^2\left[1+\dfrac{\zeta_0}{\zeta_0^*(Y;\l)}
				S^{-1,0}(x,(\l\zeta_0^*(Y;\l))^\frac{1}{m})\right]^2
		\sim\norm{x}^2,\;
		\| M \| \sim \norm{x}.
		\end{aligned}
	\end{equation}
	By \eqref{disSG}, \eqref{eq:zeta0star}, and $\zeta_0^*=G(\zeta_0^*;Y;\l)$, 
	$(X,Y)\in S_X\times S_Y\supseteq\supp\,V_2(.;\l)$, we also find
	\[
		|X_0^*(Y;\l)-X_0(Y)|=|\zeta_0^*(Y;\l)-\zeta_0(Y)|=
		|\zeta_0\,S^{-1,0}(x,(\lambda\zeta_0^*(Y;\l))^\frac{1}{m})|\le
		\frac{A\varepsilon}{2}\norm{x}^{-1},
	\]
	uniformly with respect to $\l\ge\l_0$. The proof is complete.
\end{Pf}
\begin{remark}
	The choice of $k_2$ depends only on the properties of $q$ and on the
	values of $k_1$ and $\varepsilon$, that is: we first fix $k_1>4AC>2AC>2$ and
	$\varepsilon\in\left(0,\dfrac{1}{2}\right)$, then $T>0$
	small enough as explained at the beginning of the proof of Proposition \ref{prop:i1}, 
	then $k_2>1$ as explained in the proofs of Propositions \ref{prop:j1} and \ref{prop:statF2},
	then, finally, $\l_0>2k_1\norm{2k_2}^m$.
\end{remark}
\vspace{5mm}

\noindent
The next Lemma says that the presence in the amplitude of factors which vanish
at $X=X_0^*$ implies the gain of negative powers of $\l$:

\begin{lemma}
	\label{lemma:intbyparts}
	Assume $\alpha\in\Z_+^2$, $|\alpha|>0$,
	\begin{equation}
		\label{eq:xpower}
		\begin{aligned}
		W\!=\!W(X,Y;\l)\prec 
		V_2(X,Y;\l)\,
		&t^{\alpha_1}\left[W_{\alpha_1+\alpha_2}(X,Y;\l)(\zeta-\zeta_0^*(Y;\l))^{\alpha_1+\alpha_2}\right]
		\\
		&\mbox{or}
		\\
		W\!=\!W(X,Y;\l)\prec 
		V_2(X,Y;\l)\,
		&t^{\alpha_1+\alpha_2}\left[W_{\alpha_2}(X,Y;\l)(\zeta-\zeta_0^*(Y;\l))^{\alpha_2}\right],
		\end{aligned}
	\end{equation}
	$W$ is smooth, $W_{k}(X,Y;\l)\prec\norm{x}^{k}$, $k\in\Z_+$,
	and has a $SG$-behaviour as the factors appearing in the expression of $V_2$.
	Then
	\begin{equation}
		\label{eq:intparts}
		\int e^{i\l F_2(X,Y;\l)} W(X,Y;\l)\,dX=\l^{-|\alpha|}\int e^{i\l F_2(X,Y;\l)} \widetilde{W}(X,Y;\l)\,dX,
	\end{equation}
	where $\widetilde{W}$ has the same $SG$-behaviour, 
	support and $x$-order of $V_2$, including the powers of $\zeta$.
\end{lemma}
\noindent
\begin{Pf}{.} By arguments similar to those used in the proof of Proposition \ref{prop:j1},
on $\supp\,W$
\[
	\partial_\zeta F_2(X,Y;\l)\succ \norm{x}|t|, 
	\quad
	\partial_t F_2(X,Y;\l) \succ \norm{x}|\zeta-\zeta_0^*(Y;\l)|.
\]
Assume that the first condition in \eqref{eq:xpower} holds.
Under the hypotheses, if $\alpha_1>0$, we can first insert 
$e^{i\l F_2(X,Y;\l)}=L_\zeta^{\alpha_1}e^{i\l F_2(X,Y;\l)}$ in the left hand side of \eqref{eq:intparts}, where
$L_\zeta=\dfrac{D_\zeta}{\l \, \partial_\zeta F_2(X,Y;\l)}$, and integrate by parts $\alpha_1$ times.
Similarly, if $\alpha_2>0$, we subsequently use $e^{i\l F_2(X,Y;\l)}=L_t^{\alpha_2}e^{i\l F_2(X,Y;\l)}$,
$L_t=\dfrac{D_t}{\l \, \partial_t F_2(X,Y;\l)}$, and integrate by parts $\alpha_2$ times.
The assertion then follows, remembering that $\zeta$-derivatives of $W$ produce either an
additional $\zeta^{-1}$ factor or a lowering of the exponent of $\zeta-\zeta_0^*$, and 
that $\zeta,\zeta_0^*\sim\norm{x}^{-1}$ on $\supp\,W$. The proof in the case that the second condition
in \eqref{eq:xpower} holds is the same, using first $L_\zeta$ and then $L_t$.
\end{Pf}

\vspace{5mm}

\noindent
\begin{Pf}{ of Proposition \ref{prop:i2}.}
	Define, 
	\[
		\cQ=\cQ(X,Y;\l)=\langle M(X-X_0^*(Y;\l),(X-X_0^*(Y;\l)\rangle,
	\]
	and, for $s\in[0,1]$, 
	\[
		\cF_s(X,Y;\l)=\cQ(X,Y;\l)+s \cG(X,Y;\l),
	\]
	\[
		\cG(X,Y;\l)=F_2(X,Y;\l)-\cQ(X,Y;\l).
	\]
	Remembering that $F_2(X_0^*(Y),Y;\l)=0$, $F_{2,X}^\prime(X_0^*(Y),Y;\l)=0$,
	$\cQ$ is the Taylor polynomial of degree two of $F_2$ at $X=X_0^*$, so that
	$\cG$ vanishes of order $3$ at $X=X_0^*$. Obviously, $\cF_0(X,Y;\l)=\cQ(X,Y;\l)$ and 
	$\cF_1(X,Y;\l)=F_2(X,Y;\l)$. Write
	\[
		\cJ_\tau(s)=\int e^{i\l\cF_s(X,Y;\tau^{-1})} V_2(X,Y;\tau^{-1})\,dX,
	\]
	$\tau\in(0,\l_0^{-1}]$, 
	and consider the Taylor expansion of $\cJ_\tau(s)$ of order $2\cN-1$, $\cN>1$, so that
	\[
		\left|\cJ_\tau(1)-
		\sum_{k=0}^{2\cN-1}\frac{\cJ_\tau^{(k)}(0)}{k!}\right|\le \sup_{0<s<1}\frac{|\cJ_\tau^{(2\cN)}(s)|}{(2\cN)!}.
	\]
	Since
	\[
		\cJ_\tau^{(2\cN)}(s)=(i\l)^{2\cN}\int e^{i\l\cF_s(X,Y;\tau^{-1})} \cG(X,Y;\tau^{-1})^{2\cN} V_2(X,Y;\tau^{-1})\,dX,
	\]
	Remark \ref{rem:esti2}
	and Lemma \ref{lemma:intbyparts} imply that $|\cJ_\tau^{(2\cN)}(s)|\prec\l^{-\cN}\norm{x}^{-\frac{n}{m}}$,
	$\tau\in(0,\l_0^{-1}]$, $s\in[0,1]$: indeed, it is
	easy to see, by direct computation, that $\cG$ can be bounded by linear combinations
	of expressions of the form
	\[
		t^{3}, t^2\left[W_{1}(X,Y;\tau)(\zeta-\zeta_0^*(Y;\tau))\right],
		t\left[W_{2}(X,Y;\tau)(\zeta-\zeta_0^*(Y;\tau))^2\right],
	\]
	\[
		t\left[W_{3}(X,Y;\l)(\zeta-\zeta_0^*(Y;\tau))^3\right],
	\]
	with $W_k$, $k\in\Z_+$, having the required properties.
	Then, the bound of $\cG^{2\cN}$ will always contain a term of the type
	$t^{3\cN}\left[W_{3\cN}(X,Y;\l)(\zeta-\zeta_0^*(Y;\l))^{3\cN}\right]$, which corresponds
	to the (minimun) value $|\alpha|=3\cN$ in \eqref{eq:xpower}.
	
	Each term $\cJ_\tau^{(k)}(0)$, $k=0, \dots, 2\cN-1$, has the quadratic phase function $\cQ$,
	which of course also satisfies
	\[
		\partial_\zeta\cQ(X,Y;\tau^{-1})\succ\norm{x}|t|,
		\quad
		\partial_t\cQ(X,Y;\tau^{-1})\succ\norm{x}|\zeta-\zeta_0^*(Y;\tau^{-1})|.
	\]
	Then, denoting by $\Gamma$ the Taylor expansion of $\cG$ at $X_0^*$ of order $3\cN$,
	we observe that $\cG^k-\Gamma^k$ can be bounded by polynomial expressions
	in $X-X_0^*$ of the kind appearing in the right hand side of \eqref{eq:xpower}, with $|\alpha|=\cN+k$
	(cfr. the proof of Theorem 7.7.5 in \cite{Ho1}). Setting
	\[
		\cT_\tau^{k}=\int e^{i\l\cQ(X,Y;\tau^{-1})}
		(i\l\Gamma(X,Y;\tau^{-1}))^kV_2(X,Y;\tau^{-1})\,dX,
	\]
	Lemma \ref{lemma:intbyparts} implies
	\[
		\cJ_\tau^{(k)}(0)- \cT_\tau^{k} \prec \l^{-\cN}\norm{x}^{-\frac{n}{m}}.
	\]
	
	We now apply the Stationary Phase Method to $\cT_\tau^k$ and prove that
	\begin{equation}
		\label{eq:statj2}
		\cJ_\tau(1)\sim\sum_{j=0}^{+\infty} d_j(Y;\tau) \, \l^{-1-j},
	\end{equation}
	which is a consequence of
	\begin{equation}
		\label{eq:statexp}
		\begin{aligned}
		\cT_\tau^k\sim\l^{-1}&\det(M/2\pi i)^{-\frac{1}{2}}
		\sum_{l} L_{l,k,Y,\tau}V_2,
	\\
		L_{l,k,Y,\tau}V_2=\sum_{l} (2i\l)^{-l}
		&\langle M^{-1}D_X,D_X\rangle^l
		[(i\l\Gamma)^k V_2](X_0^*(Y;\tau^{-1}),Y;\tau^{-1})/l!\,,
		\end{aligned}
	\end{equation}
	with $M$ evaluated with $\tau^{-1}$ in place of $\l$.
	Recalling \eqref{eq:M}, it follows that the inverse matrix $M^{-1}$ satisfies, on $\supp\,V_2$,
	\[
		M^{-1}=
		\left(
			\begin{array}{cr}
				0
				&
				\dfrac{1}{M_{12}}
				\\
				\dfrac{1}{M_{12}}
				&
				-\dfrac{M_{11}}{M_{12}^2}\rule{0mm}{6mm}
			\end{array}
		\right),
		\quad
		\dfrac{1}{M_{12}}\prec\norm{x}^{-1}, \dfrac{M_{11}}{M_{12}^2}\prec\varepsilon\norm{x}^{-2},
		\quad\|M^{-1}\|\sim\norm{x}^{-1},
	\]
	in view of the ellipticity of the involved symbols. Then, the operators $L_{j,k,Y,\tau}$, $j,k\in\Z_+$,
	do not increase the $x$-order of the resulting function with respect to that of their arguments,
	$(i\l\Gamma)^k V_2$, which is the same of $V_2$, uniformly with respect to $\tau$. The proof
	of \eqref{eq:statexp} then follows
	by Theorem 7.6.1, the proof of Lemma 7.7.3 
	and formula (7.6.7) in \cite{Ho1}, see also \cite{HeRo:1,HeRo81}. Indeed, by the mentioned
	results, 
	\begin{align*}
		\cJ_\tau^k&  -  
		\l^{-1}\det(M/2\pi i)^{-\frac{1}{2}}\sum_{l\le k+\cN}
		L_{l,k,Y,\tau}V_2
		\\
		&=\left[\cJ_\tau^k - 
		\l^{-1}\det(M/2\pi i)^{-\frac{1}{2}}\sum_{l\le k+\cN+1}
		L_{l,k,Y,\tau}V_2\right]
		+\l^{-1}\det(M/2\pi i)^{-\frac{1}{2}}%\sum_{k+\cN+1\le l\le k+\cN+2}
		L_{k+\cN+1,k,Y,\tau}V_2
		\\
		&\prec \l^{-\cN-3}\norm{x}^{-1}\sum_{|\beta|\le2}\| 
		D^\beta_X \langle M^{-1}D_X,D_X\rangle^{k+\cN+3} [\cG^k V_2](X,Y;\tau^{-1})
		\|_{L^2(\R^2_X)}
		\\
		&+\l^{-\cN-2}\norm{x}^{-1}|L_{k+\cN+1,k,Y,\tau}V_2(X_0^*(Y;\tau^{-1}),Y;\tau^{-1})|
		%+\l^{-1}L_{k+\cN+2,k,Y,\tau}V_2(X_0^*(Y;\tau^{-1}),Y;\tau^{-1})|
		\\
		&\prec \l^{-\cN-3}\norm{x}^{-1}\left[ 
			\int_{c^{-1}\norm{x}^{-1}}^{c\norm{x}^{-1}} \zeta^{2\left(\frac{n}{m}-3\right)}\,d\zeta
		\right]^\frac{1}{2}+
		\l^{-\cN-2}\norm{x}^{-\frac{n}{m}}
		\\
		&\prec\l^{-\cN-3}\norm{x}^{-\frac{n}{m}+\frac{3}{2}}+\l^{-\cN-2}\norm{x}^{-\frac{n}{m}}
		\\
		&\prec\l^{-\cN-1-\frac{1}{2}}\norm{x}^{-\frac{n}{m}}, \quad\l\to+\infty,
	\end{align*}
	since $\norm{x}\prec\l$ on $\supp\,V_2$. It is then enough to sum all the expansions
	of $\dfrac{\cJ_\tau^k}{k!}$, $k=0,\dots,2\cN-1$, and sort the terms
	by decreasing exponents of $\l$ (as in the proof of Theorem 7.7.5 in \cite{Ho1})
	to obtain \eqref{eq:statj2} with the usual expression
	\[
		\widetilde{d}_j(Y;\tau)=\det(M/(2\pi i))^{-\frac{1}{2}}\sum_{k-l=j}\sum_{2k\ge3l}i^{-j}2^{-k}
		\langle M^{-1}D_X,D_X\rangle^{k} [(i\Gamma)^l V_2](X_0^*(Y;\tau^{-1}),Y;\tau^{-1}),
	\]
	so that, in particular,
	\[
		\widetilde{d}_j(Y;\tau)\prec\norm{x}^{-\frac{n}{m}},
	\]
	for any $j\in\Z_+$, $\tau\in(0,\l_0^{-1}]$. We can then integrate $\cJ_\tau(1)$ and its 
	asymptotic expansions with respect to $Y\in\Sn\times\{x\in\R^n\colon\norm{x}\le\varkappa\l\}$ and find
	\begin{equation}
		\label{eq:proof2nd-1}
		J_2(\l)=\int_{\norm{x}\le\varkappa\l}\int_{\Sn} 
		\cJ_{\l^{-1}}(1)\,dY\sim \sum_j \l^{-1-j}
		\int_{\norm{x}\le\varkappa\l}\int_{\Sn} \widetilde{d}_j(Y;\l^{-1})\,dY, \quad\l\to+\infty,
	\end{equation}
	recalling that 
	$\psi(0)=1$ and $a(0,x,\xi)=1$, for all $x,\xi\in\R^n$.
	Moreover, for $\zeta=\zeta_0^*(Y;\l)$, the factors $H_1$, $H_2$, and $H_3$ are 
	identically equal to $1$ (see the Appendix).
	Then, the coefficient of the leading term in \eqref{eq:proof2nd-1}
	is given by 
	\begin{align*}
		&\int\widetilde{d}_0(Y;\l^{-1})\,dY=\int_{\norm{x}\le\varkappa\l}\int_{\Sn}
		\det(M/(2\pi i))^{-\frac{1}{2}} \, V_2(X_0^*(Y;\l),Y;\l)\,dY
		\\
		&=2\pi\int_{\norm{x}\le\varkappa\l}\int_{\Sn}    
		H_1\!\left(\frac{\langle x \rangle \langle (\l\zeta_0^*(\varsigma,x;\l))^\frac{1}{m} \rangle^m}{\l}\right)
		H_2((\l\zeta_0^*(\varsigma,x;\l))^\frac{1}{m}) \, H_3\left(\dfrac{\zeta_0^*(\varsigma,x;\l)}{\zeta_0(\varsigma,x)}-1\right) \cdot
		\\
		&\hspace*{27mm}\cdot |\det(M)|^{-\frac{1}{2}}\,\zeta_0^*(\varsigma,x;\l)^{\frac{n}{m}-1}
		\,d\varsigma dx
		\\
		&=2\pi\int_{\norm{x}\le\varkappa\l}\int_{\Sn}    
		\;\;|\det(M)|^{-\frac{1}{2}}\,
		\zeta_0^*(\varsigma,x;\l)^{\frac{n}{m}-1}
		\,d\varsigma dx,
	\end{align*}
	with $M$ evaluated in $\zeta=\zeta_0^*$. We say that 
	\begin{align*}
		\int\widetilde{d}_0(Y;\l^{-1})\,dY &= 
		2\pi\int_{\R^n}\int_{\Sn}\frac{1}{q_\psi(x,\varsigma)^\frac{n}{m}}\,d\varsigma dx
		+O(\l^{\max\{-\frac{1}{m}, n-\frac{n}{m}, -1\}})
		\\
		&=2\pi d_0
		+O(\l^{\max\{-\frac{1}{m}, n-\frac{n}{m}, -1\}}),
		\quad\l\to+\infty.
	\end{align*}
	To confirm this, first note that $\zeta_0^*(Y;\l)\to\zeta_0(Y)$, $\l\to+\infty$, for any $(Y;\l)$ belonging to
	the support of the integrand, see the Appendix. Moreover, 
	the integrand is uniformly bounded by the summable function $\norm{x}^{-\frac{n}{m}}$, and its support
	is included in the set $S$.
	Then, recalling \eqref{eq:M} and setting $\widetilde{H}=|\zeta_0^2\det(M)|^{-\frac{1}{2}}$,
	\begin{align*}
		R&=\int_{\norm{x}\le\varkappa\l}\int_{\Sn}   
		|\det(M)|^{-\frac{1}{2}}\,
		\zeta_0^*(Y;\l)^{\frac{n}{m}-1}\,dY
		-\int_{\R^n}\int_{\Sn}   
		\zeta_0(Y)^{\frac{n}{m}-1}\,dY
		\\
		&=\int_{\norm{x}\le\varkappa\l}\int_{\Sn}    
		 \zeta_0 \left[\widetilde{H}\,
		(\zeta_0^*)^{\frac{n}{m}-1}
		-\zeta_0^{\frac{n}{m}-1}\right]\,d\varsigma dx
		-\int_{\norm{x}\ge\varkappa\l}\int_{\Sn}    
		\zeta_0^\frac{n}{m}\,d\varsigma dx.
	\end{align*}
	The second integral is always $O(\l^{n-\frac{n}{m}})$, since $q_\psi(x,\varsigma)\sim\norm{x}$ implies
	\begin{align*}
		R_2=\int_{\norm{x}\ge\varkappa\l}\int_{\Sn}    
		\zeta_0^\frac{n}{m}\,d\varsigma dx &\sim \int_{\varkappa\l}^{+\infty}r^{n-\frac{n}{m}-1}\,dr
		=\frac{(\varkappa\l)^{n-\frac{n}{m}}}{\dfrac{n}{m}-n},
		\quad\l\to+\infty.
	\end{align*}
	The first integral can be estimated as follows. Since 
	\[
		\zeta_0^*-\zeta_0=\zeta_0 S^{-1,0}(x,(\l\zeta_0^*)^\frac{1}{m})=\zeta_0 O((\l\zeta_0^*)^{-\frac{1}{m}}),
	\]
	by the properties of $\zeta_0^*$ (see the Appendix) we find
	\[
		\left(\frac{\zeta_0^*}{\zeta_0}\right)^{\frac{n}{m}-1}-1 = (1+O((\l\zeta_0^*)^{-\frac{1}{m}}))^{\frac{n}{m}-1}-1=
		O((\l\zeta_0^*)^{-\frac{1}{m}}))=O(\norm{x}^\frac{1}{m}\l^{-\frac{1}{m}}),
	\]
	since $S^{-1,0}(x,(\l\zeta_0^*)^\frac{1}{m})<<1$. By \eqref{eq:M}, we similarly have
	$\widetilde{H}=1+O(\norm{x}^\frac{1}{m}\l^{-\frac{1}{m}})$, so that
	\begin{align*}
		R_1=\int_{\norm{x}\le\varkappa\l}\int_{\Sn}    
		 \zeta_0 &\left[\widetilde{H}\,
		(\zeta_0^*)^{\frac{n}{m}-1}
		-\zeta_0^{\frac{n}{m}-1}\right]\,d\varsigma dx
		\\
		&=
		\int_{\norm{x}\le\varkappa\l}\int_{\Sn}    
		 \zeta_0^\frac{n}{m} \left[\widetilde{H}\,
		\left(\frac{\zeta_0^*}{\zeta_0}\right)^{\frac{n}{m}-1}
		-1\right]\,d\varsigma dx
		\\
		&\prec
		\l^{-\frac{1}{m}}
		\int_{\norm{x}\le\varkappa\l}   
		 \norm{x}^{-\frac{n-1}{m}}\, dx.
	\end{align*}
	If $n>\dfrac{1}{1-m}\Leftrightarrow n-1-\dfrac{n-1}{m}<-1$, $n\in\N$, $m\in(0,1)$,
	the integral in $R_1$ is convergent for $\l\to+\infty$ and $R_1=O(\l^{-\frac{1}{m}})$.
	In this case, $R_1$ contributes an $O(\l^{\frac{n}{m}-1-\frac{1}{m}})$ term
	to the expansion of $I_2(\l)$, which is of lower order than the $O(\l^{\frac{n}{m}-2})$ term,
	which is one of the remainders appearing in \eqref{eq:i2weyl}.
	On the other hand, if $n<\dfrac{1}{1-m}$, the integral in $R_1$ is divergent, and $R_1$ itself
	is $O(\l^{n-\frac{n}{m}})$, since, trivially,
	\[
		\lim_{\l\to+\infty}\frac{\displaystyle\l^{-\frac{1}{m}}
		  \int_0^{\varkappa\l}\dfrac{r^{n-1}}{(1+r^2)^\frac{n-1}{2m}}dr}{\l^{n-\frac{n}{m}}}
		=\lim_{\l\to+\infty}\frac{\displaystyle
		  \int_0^{\varkappa\l}\dfrac{r^{n-1}}{(1+r^2)^\frac{n-1}{2m}}dr}{\l^{n-\frac{n-1}{m}}}
		=\frac{\varkappa^{n-1-\frac{n-1}{m}}}{n-\dfrac{n-1}{m}}.
	\]
	Finally, if $n=\dfrac{1}{1-m}$, $R_1$ is $O(\l^{-\frac{1}{m}}\ln\l)$, by
	\[
		\lim_{\l\to+\infty}\frac{\displaystyle
		  \int_0^{\varkappa\l}\dfrac{r^{\frac{1}{1-m}-1}}{(1+r^2)^\frac{1}{2(1-m)}}dr
		  }{\ln\l}
		=\varkappa^{-1},
	\]
	and again contributes a term of lower order than the remainder $O(\l^{\frac{n}{m}-2})$.
	Similar conclusions can be obtained for the subsequent terms of the expansion of $J_2(\l)$.
	The proof is complete, combining the contributions of the remainders
	like $R$ with the other terms in
	the expansion of $J_2(\l)$, and remembering that
	\[
		I_2(\l)=\frac{n}{m}\frac{\l^\frac{n}{m}}{(2\pi)^n}J_2(\l)+O(|\l|^{-\infty})
		=\frac{n}{m}\l^\frac{n}{m}
		 \sum_{j=0}^{+\infty} (d_j\l^{-1-j}+O(\l^{n-\frac{n}{m}-1-j}))+O(|\l|^{-\infty}).
	\]
\end{Pf}

\begin{remark}
	The same conclusions concerning the behaviour of $R_1$ in the final step of 
	the proof of Proposition \ref{prop:i2} could have been obtained studying the 
	Taylor expansion of the extension of $\zeta_0^*(Y;\tau^{-1})$, $\tau=\l^{-1}$, 
	to the interval $[0,\l_0^{-1}]_\tau$, similarly to \cite{HeRo81}. 
\end{remark}
\vspace{5mm}
\noindent
\begin{Pf}{ of Theorem \ref{thm:main}.}
	The statement for $\mu>m$ follows by the arguments in Section \ref{sec:N_a}
	and Propositions \ref{prop:h1}, \ref{prop:i1}, \ref{prop:i2},
	summing up the contribution of the local symbol on the exit chart
	to the contributions of the remaining local symbols, which
	gives the desired multiple of the integral of $q_\psi^{-\frac{n}{m}}$ on the cosphere bundle
	as coefficient of the leading term $\l^\frac{n}{m}$. The remainder has then order
	equal to the maximum between $\dfrac{n}{m}-1$ and $n$, as claimed.
	The proof for $\mu<m$ is the same, by exchanging step by step the
	role of $x$ and $\xi$.
\end{Pf}

\section*{Appendix}
\label{sec:app}
\setcounter{equation}{0}

For the sake of completeness, here we illustrate some details of the proof of Theorem \ref{thm:main},
which we skipped in the previous Sections. They concern, in particular, formula \sref{eq:tr3}, which 
expresses the relation between 
$\displaystyle\sum_{j} \widehat{\psi}(\lambda-\eta_j)$ and the oscillatory integrals examined in Section \ref{sec:stat}. We mainly focus on the aspects which are specific for the manifolds with ends.
We also show more precisely how the constants $k_1,k_2,\l$ are involved in the solution
of equation \eqref{eq:zeta0star} via the Fixed Point Theorem, 
completing the proof of Proposition \ref{prop:statF2}.

\renewcommand{\thesection}{A}

\subsection{Solution of Cauchy problems and $SG$ Fourier Integral Operators}
\label{sec:a1.2}
$ $\\
Using the so-called ``geometric optics method'', specialised to che 
pseudodifferential calculus we use (see 
\cite{Coriasco:998.1, Coriasco:998.2,Coriasco-Maniccia:003,Coriasco-Panarese:001,Coriasco:998.3} and
\cite{Nicola-Rodino:003}), the Cauchy Problem \sref{eq:CP} can be 
solved modulo $\cS(M)$ by means of an operator family $V(t)$, defined 
for $t$ in a suitable interval $(-T, T)$, 
$T > 0$: $V(t)$ induces continuous maps
\begin{eqnarray*}
	& & V\;\colon\;\cS(M) \rightarrow C^\infty((-T,T), \cS(M)),
	\\
	& & V\;\colon\;\cS^\prime(M) \rightarrow C^\infty((-T,T), \cS^\prime(M))
\end{eqnarray*}
and
\begin{eqnarray}
	\label{eq:CP1}
	& &
	(D_t-Q)\circ V =: R \in C^\infty((-T,T), L^{-\infty,-\infty}(M)),
	\\
	\label{eq:CP2}
	& &
	Vu_{| t=0} = u, \hspace{1cm} \forall u \in \cS^\prime(M).
\end{eqnarray}
First of all, we recall that the
partition of unity $\{\theta_k\}$ and the family of functions
$\{\chi_k\}$ of Definition \ref{def:mwce} can be chosen so that 
$(\theta_k)_*$ and $(\chi_k)_*$ are $SG$-symbols of
 order $(0,0)$ on $U_k$, extendable to symbols of the same class
 defined on $\R^n$ (see \cite{SC87}).

\begin{remark}
\label{rem:varphi}
\begin{enumerate}
\item The complete symbol of $Q$ depends, in general, on the choice 
of the admissible atlas, 
          of $\{ \theta_k \}$ and of $\{ \chi_k \}$. Anyway, if $\{ 
\tilde{q}_k \}$ is another complete symbol of
          $Q$, $\kappa(x) (q_k(x,\xi) - \tilde{q}_k(x,\xi)) \in 
\cS(\varphi_k(\Omega_k \cap \widetilde{\Omega}_k))$ for an
          admissible cut-off function $\kappa$ supported in  
$\varphi_k(\Omega_k \cap \widetilde{\Omega}_k)$.
\item The solution of \sref{eq:CP} in the $SG$-classical case and 
the properties of $\varphi_k$ and $a_k$ 
	in \sref{eq:ph-amp} were investigated in 
\cite{Coriasco-Panarese:001} (see also \cite{Nicola-Rodino:003}, Section 4). 
In particular, it turns out that 
	$\varphi_k \in C^\infty((-T_k,T_k), S_\cl^{1,1})$,
	$T_k > 0$.
	According to \cite{Coriasco:998.2}, page 101, for every $SG$ phase 
functions $\varphi$  
	of the type involved in the definition of $V(t)$ we also have, for 
all $x\in\R^n$:
	\begin{eqnarray*}
		|\nabla_\xi\varphi(t;x,\xi) - x| 
		& = & |\nabla_\xi\varphi(t;x,\xi) - \nabla_\xi\varphi(0;x,\xi)|
		   =     \left| 
		   	        \int_0^t \nabla_\xi\dot{\varphi}(t;x,\xi)\, dt
		          \right|
		\\
		& = & \left| 
		   	        \int_0^t \nabla_\xi(q(x, d_x\varphi(t;x,\xi))) \, dt 
			 \right| \le C|t|\norm{x},
	\end{eqnarray*}
	with a constant $C>0$ not depending on $t,x,\xi$. The function
	$\Phi_{t,\xi}(x):= \nabla_\xi\varphi(t,x,\xi)$ turns out to be a 
($\!SG$-)diffeomorphism, smoothly 
	depending on the parameters $t$ and $\xi$ (see 
\cite{Coriasco:998.1}).
\end{enumerate}
\end{remark}
Before proving Theorem \ref{thm:globalV}, we state a technical Lemma,
whose proof is immediate and henceforth omitted. %
\begin{lemma}
	\label{lemma:supp}
	Let $U \subset \R^n$ be an open set and define 
	$\displaystyle U_{\delta} := \bigcup_{x\in U}
	B(x,\delta\norm{x})$ for arbitrary $\delta>0$.
	Assume 
	$\theta, \chi \in C^\infty(\R^n)$ such that
	$\displaystyle \mathrm{supp}\, \theta \subset 
U_{\frac{\delta}{5}}$,
	$\displaystyle \mathrm{supp}\, \chi \subset U_{\delta}$ and
	$\displaystyle \chi|_{U_\frac{\delta}{2}} \equiv 1$. 
	Then, for any diffeomorphism $\Phi_{t, \xi}$, smoothly depending on 
$t\in (-T,T)$, $\xi\in\R^n$, and such that
	$\displaystyle \forall t,x,\xi \;  |\Phi_{t,\xi}(x) - x| \le 
C|t|\norm{x}$ with a  constant $C>0$ independent 
	of $t,x,\xi$,
	\[
	|t| < \frac{\delta}{4C} \Rightarrow (1-\chi(x)) \, 
(\partial^\alpha\theta)(\Phi_{t,\xi}(x)) = 0
	\] 
	for any multiindex $\alpha$ and $x,\xi\in\R^n$.
\end{lemma}

\noindent
We remark that, since  a manifold with ends is, in particular, a $SG$-manifold, 
the charts $(\Omega_k,\psi_k)$ and the functions $\{\theta_k\}$, $\{\chi_k\}$,
can be chosen such that
\begin{itemize}
	\item for a fixed $\delta>0$, each coordinate open set $U_k=\psi_k(\Omega_k)$,
	$k=1,\dots,N$,
	contains an open subset $W_k$ such that
	$\displaystyle \bigcup_{x\in W_k}B(x,\delta\norm{x})\subseteq U_k$;
	\item the supports of $\theta_k$ and $\chi_k$, $k=1,\dots,N$, satisfies
	hypotheses as the supports of $\theta$ and $\chi$ in Lemma \ref{lemma:supp} (see, e.g., Section 3 of 
	\cite{SC87} for the construction of functions with the required properties). 
\end{itemize}
\noindent
In fact, this is relevant only for $k=N$.

\medskip

\noindent
\begin{Pf}{ of Theorem \ref{thm:globalV}.}
	We will write $R \equiv S$ when $R - S \in L^{-\infty,-\infty}(M)$ 
and $\chi_k \vartriangleleft 
	\widetilde{\chi}_k$ when the functions $\chi_k, 
\widetilde{\chi}_k$ are smooth, 
	non-negative, supported in $\Omega_k$, satisfy $\chi_k \, 
\widetilde{\chi}_k = \chi_k$ 
	and $(\chi_k)_*, (\widetilde{\chi}_k)_* $ are $SG$-symbols
	of order $(0,0)$ on $U_k=\psi_k(\Omega_k)$. Obviously, $R \in 
L^{-\infty,-\infty}(M)$ implies $R\,V \in 
	C^\infty((-T,T), L^{-\infty,-\infty}(M))$. To 
simplify notation, in the computations 
	below we will not distinguish between the functions $\chi_k$, 
$\theta_k$, etc., and their local 
	representations.
	
	$V(t)$ obviously satisfies \eqref{eq:CP2}. To prove \eqref{eq:CP1},
	choose functions $\zeta_k, \upsilon_k$ supported in $\Omega_k$ such that
	$\theta_k \vartriangleleft \zeta_k \vartriangleleft \chi_k  \vartriangleleft \upsilon_k$. Then
	$\displaystyle Q \equiv \sum_{k=1}^N \theta_k Q_k \chi_k$ and, for 
all $k = 1, \dots, N$,
	$Q \chi_k \equiv \upsilon_k Q_k \chi_k$  
	(see \cite{CO}, Section 4.4; cfr. also \cite{Kumano-go:1}), so 
that
	\begin{eqnarray}
		\nonumber
		Q \, V(t)  & = & \sum_{k=1}^N  Q \, \chi_k \, V_k(t) \, \theta_k \equiv
		                                  \sum_{k=1}^N  \upsilon_k \, Q_k \, 
\chi_k \, V_k(t) \theta_k
		                                  \\
		\nonumber
		               & = & \sum_{k=1}^N  (\upsilon_k \, [Q_k, \chi_k]  \, 
V_k(t) \, \theta_k 
		                                   + \chi_k \, Q_k \, V_k(t) \, 
\theta_k)
		                                  \\
		\label{eq:dtvt}
		               & \equiv  & \sum_{k=1}^N  ( \upsilon_k \, [Q_k, \chi_k]  
\, 
\zeta_k       \, V_k(t) \, \theta_k 
		                                   + \upsilon_k \, [Q_k, \chi_k]  
\, (1-\zeta_k) \, V_k(t) \, \theta_k)
		                                   + D_t V(t) \equiv D_t V(t).
	\end{eqnarray}
	That the first term in the sum \eqref{eq:dtvt} is smoothing comes 
from the $SG$ symbolic calculus
	in $\R^n$ and the observations above,
	since $\mathrm{sym}\,([Q_k, \chi_k]  \, \zeta_k)$ $\sim 0$. The 
same property holds for
	each $k$ in the second term, provided $t\in I_{T_k}$, $T_k > 0$ 
	small enough. In fact, by Theorems 7 and 8 of \cite{Coriasco:998.1}, 
$(1-\zeta_k)\,V_k(t)\, 
	\theta_k$ is a $SG$ 
	FIO with the same phase function $\varphi_k$ and amplitude $w_k$ 
such that
	\[
	w_k(x,\xi) \sim \sum_\alpha 
	\frac{(1-\zeta_k(x)) \, 
(\partial^\alpha\theta_k)(\nabla_\xi\varphi_k(t;\xi,x))}{\alpha!}
	                               		b_{j\alpha}(t; x,\xi),
	\]
	with suitable $SG$-symbols $b_{j\alpha}$ defined in terms of 
$\varphi_k$ and $a_k$. By 
	Remark \ref{rem:varphi} and Lemma
	 \ref{lemma:supp}, $w_k \sim 0$ for $|t|$ small enough.
	 The proof that $V(t)$ satisfies \eqref{eq:CP1} is completed once we 
set 
	 $T = \min \{ T_1, \dots, T_N\}$. The 
last part of the Theorem can be
	 proved as in \cite{GriSjo:994}, Proposition 12.3, since, setting
	 $W(t) := U(-t) \, V(t)$, it is easy to see $D_t W(t) \equiv 0$, so 
that 
	 $W(0) = I  \Rightarrow
	 W(t) \equiv I \Rightarrow V(t) \equiv U(t)$, with smooth dependence 
on $t$, as claimed.
\end{Pf}

\subsection{Trace formula and asymptotics for $A\in EL_\cl^{r,1}(M)$}
\label{sec:a1.3}
$ $\\
\begin{Pf}{ of Lemma \ref{lemma:kernconv}.}
	Consider first the finite sum
	\[
	k_J(x,y) = \sum_{j=1}^J \widehat{\psi}(-\eta_k) e_k(x)  \overline{e_k(y)}
	\]
	and reduce to the local situation (cfr. Schrohe \cite{SC87}), via the $SG$-compatible partition 	
	of unity $\{\theta_l\}$ subordinate to the atlas $\cA$, by
	\[
	k_J(x,y) = \sum_{r,s=1}^N 
	                    \sum_{j=1}^J \widehat{\psi}(-\eta_k) \, (\theta_r e_k)(x) \, \overline{(e_k \theta_s)(y)}
		       = \sum_{r,s=1}^N k^{rs}_R(x,y).
	\]
	Then, by $e_k \in \cS(M)$ and the fact that $(\theta_r)_* = \theta_r\circ \psi_r^{-1}$ is 
	supported and at most of polynomial growth in $U_r$, it turns out that we can extend 
	$(\theta_r e_k)_*$ and $(\theta_s e_k)_*$ to elements of  $\cS(\R^n)$.
	By an argument similar to the proof of Proposition 1.10.11
	in \cite{Helffer:984.1} (or by direct estimates of the involved seminorms, as in \cite{GriSjo:994}), 
	$(k_J^{rs})_* \rightarrow (k^{rs})_*$ in $\cS(\R^n\times\R^n)$ when $J\to+\infty$, with 
	$(k^{rs})_*$ kernel of $(\theta_r \widehat{\psi}(-Q) \, \theta_s)_*$. This proves that 
	$\displaystyle \widehat{\psi}(-Q) = \sum_{r,s=1}^N \theta_r  \widehat{\psi}(-Q) \, \theta_s$
	is an operator with kernel 
	$\displaystyle K_\psi(x,y) = \sum_{r,s=1}^N k^{rs}(x,y) \in \cS(M\times M)$.
\end{Pf}

\noindent
The proof of Theorem \ref{thm:taub} is essentially the one in \cite{GriSjo:994},
while the proof of Lemma \ref{lemma:R} comes from \cite{Helffer:984.1}: we include
both of them here for convenience of the reader.

\noindent
\begin{Pf}{ of Theorem \ref{thm:taub}.}
Setting $\displaystyle G(\lambda)=\int_{-\infty}^\lambda \widehat{\psi}(\tau) d\tau$ and integrating \sref{eq:asint1} in $(-\infty,\lambda)$, we obtain
\begin{eqnarray}
\nonumber
\lefteqn{
\int_{-\infty}^\lambda \int \widehat{\psi}(\tau-\eta) dN_Q(\eta) d\tau =
\int \left( \int_{-\infty}^\lambda \widehat{\psi}(\tau-\eta) d\tau\right) dN_Q(\eta)
}
\\
& &
\label{eq:asint2}
\hspace{1cm}
= \int G(\lambda-\eta) dN_Q(\eta) =
	\left\{
	\begin{array}{ll}
		\displaystyle d_0\, \lambda^{\frac{n}{m}} + O(\lambda^{n^*}) &
		 \mbox{for } \lambda\to+\infty
		\\
		\displaystyle\rule{0mm}{7mm}O(|\lambda|^{-\infty}) & \mbox{for } \lambda\to-\infty.
	\end{array}
	\right. 
\end{eqnarray}
Now, observe that
$$
\int G(\lambda - \eta) dN_Q(\eta) = \sum_{j=1}^\infty G(\lambda-\eta_j) =
\sum_{j=1}^\infty \int_{-\infty}^{\lambda-\eta_j} \widehat{\psi}(\tau) d\tau
= \sum_{j=1}^\infty \int H(\lambda-\eta_j-\tau) \widehat{\psi}(\tau) d\tau,
$$
where $H(\tau)$ is the Heaviside function. Bringing the series under the integral sign, we can write
\begin{equation}
\label{eq:asint3}
\begin{array}{l}
\hspace{-0.4cm}
\displaystyle
\int G(\lambda-\eta) dN_Q(\eta) = 
\int \sum_{j=1}^\infty H(\lambda-\eta_j-\tau) \widehat{\psi}(\tau) d\tau
= \int N_Q(\lambda-\tau) \widehat{\psi}(\tau) d\tau
\\
\hspace{0.7cm}
\displaystyle\rule{0mm}{7mm}
= N_Q(\lambda) \int \widehat{\psi}(\tau) d\tau + \int [N_Q(\lambda-\tau) - N_Q(\lambda)] \widehat{\psi}(\tau) d\tau
= 2\pi N_Q(\lambda) + R(\lambda),
\end{array}
\end{equation}
since $\displaystyle\int\widehat{\psi}(\tau)d\tau=2\pi\psi(0)=2\pi$.
In view of the monotonicity of $N$ and next Lemma \ref{lemma:R} (cfr. Lemma 4.2.8 
of \cite{Helffer:984.1}), for $\lambda\ge 1$
\begin{eqnarray*}
\lefteqn{
| N_Q(\lambda-\tau) - N_Q(\lambda)| \le N_Q(\lambda + |\tau|) - N_Q(\lambda - |\tau|)
= \int_{\lambda-|\tau|}^{\lambda+|\tau|} dN_Q(\eta) 
}
\\
& & \hspace{1cm}
= \int_{|\lambda-\eta| \le |\tau|} dN_Q(\eta)
\le C (1+|\tau|)^\frac{n}{m} (1 + |\lambda|)^{\frac{n}{m}-1}
\le \tilde{C} (1 + |\tau|)^\frac{n}{m} \lambda^{\frac{n}{m}-1}.
\end{eqnarray*}
We can then conclude that $R(\lambda)=O( \lambda^{\frac{n}{m}-1})$, $\lambda\ge 1$, since $\widehat{\psi}\in\cS$, and
this, together with \sref{eq:asint2} and \sref{eq:asint3}, completes the proof.
\end{Pf}

\begin{lemma}
	\label{lemma:R}
	Under the hypotheses of Theorem \ref{thm:taub},
	there exists a constant $C>0$ such that for any $K\ge 0$ and any $\lambda\in\R$
	$$
	\int_{|\lambda-\eta|\le K} dN_Q(\eta) \le C (1+K)^\frac{n}{m} (1+|\lambda|)^{\frac{n}{m}-1}
	$$
\end{lemma}
\noindent
\begin{Pf}{.}
Let $h\in(0,\widehat{\psi}(0))$ and $[-K_0,K_0]$ such that $\widehat{\psi}(t)\ge h$ for all $t\in[-K_0,K_0]$. Then, trivially,
$$
	\int_{|\lambda-\eta|\le K_0} dN_Q(\eta) \le h^{-1} \int \widehat{\psi}(\lambda-\eta) dN_Q(\eta).
$$
Let us now prove that
$$
\int \widehat{\psi}(\lambda-\eta) dN_Q(\eta) \le C_1(1+|\lambda|)^{\frac{n}{m}-1}.
$$
Indeed, this is clear for $\lambda\ge\tilde{C}>0$ and $\lambda\le-\tilde{C}$, $\tilde{C}$ suitably large, in view of hypothesis (iii). For $\lambda\in[-\tilde{C},\tilde{C}]$, choose a constant $C_1$ so large that
$\displaystyle \max_{\lambda\in[-\tilde{C},\tilde{C}]}\int \widehat{\psi}(\lambda-\eta) dN_Q(\eta) \le C_1(1+\tilde{C})^{\frac{n}{m}-1}$. This shows that, for all $\lambda\in\R$,
\begin{equation}
	\label{eq:est1}
	\int_{|\lambda-\eta|\le K_0} dN_Q(\eta) \le C_2(1+|\lambda|)^{\frac{n}{m}-1}.
\end{equation}
For arbitrary $K>0$ there exists $l\in\N$ such that $(l-1)K_0 \le K < lK_0$. We write
\begin{eqnarray*}
\lefteqn{
\int_{|\lambda-\eta|\le K} dN_Q(\eta) \le \int_{|\lambda-\eta|\le lK_0} dN_Q(\eta)
=\sum_{j=0}^{l-1} \int_{jK_0\le|\lambda-\eta|\le (j+1)K_0} dN_Q(\eta)
}
\\
& & \hspace{1.8cm}
\le \sum_{j=0}^{l-1}
\left[
\int_{\left|\lambda+jK_0+\frac{K_0}{2}-\eta\right|\le \frac{K_0}{2}} dN_Q(\eta)
+
\int_{\left|\lambda-jK_0-\frac{K_0}{2}-\eta\right|\le \frac{K_0}{2}} dN_Q(\eta)
\right].
\end{eqnarray*}
By \sref{eq:est1}, the last sum can be estimated by 
\begin{eqnarray*}
\lefteqn{
2C_2\sum_{j=0}^{l-1} \left(1+|\lambda| +\left( j+\frac{1}{2}\right)K_0  \right)^{\frac{n}{m}-1}
\le
2C_2l \left(1+|\lambda| +\dfrac{K_0}{2}+K  \right)^{\frac{n}{m}-1}
}
\\
& & \hspace{1cm}
\le
2C_2\left(1+\frac{K}{K_0}\right) \left(1+\frac{K_0}{2}+K\right)^{\frac{n}{m}-1}
\left(1+\frac{1}{1+\dfrac{K_0}{2}+K} |\lambda|\right)^{\frac{n}{m}-1}
\\
& & \hspace{1cm}
\le C (1+K)^\frac{n}{m} (1+|\lambda|)^{\frac{n}{m}-1},
\end{eqnarray*}
as claimed.
\end{Pf}

\subsection{The solution $\zeta_0^*(Y;\l)$ of the equation $\zeta=G(\zeta;Y;\l)$.}
\label{sec:a1.4}
$ $\\
We know that $A^{-1}\norm{x}^{-1}\le \zeta_0(\varsigma, x)=q_\psi(x,\varsigma)^{-1}\le A\norm{x}^{-1}$,
$Y=(\varsigma,x)\in \widetilde{S}_Y=\Sn\times\{x\in\R^n\colon\norm{x}\le\varkappa\l\}$,
and that $k_1>4AC>2AC>2$.
Moreover, $k_2>1$ is chosen so large that, in particular, on $\supp\,U_2\supset\supp\,V_2$,
the absolute value of the $\zeta$-derivative of $G$ is less than $k_0\le1$,
uniformly with respect to $Y\in \widetilde{S}_Y$, $\l\ge\l_0$, $(X,Y;\l)\in\supp\,V_2$. 
We want to show that once $k_1$ is fixed, the choice of such a suitably large $k_2>1$
allows to make  $G$ a contraction on the compact set
$I_x=\left[A^{-1}\left(1-\dfrac{\varepsilon}{2}\right)\norm{x}^{-1}, A\left(1+\dfrac{\varepsilon}{2}\right)\norm{x}^{-1}\right] \subset[c^{-1}\norm{x}^{-1}, c\norm{x}^{-1}]$, uniformly with respect to 
$(\varsigma,x)$, 
$\l\ge\l_0$, provided $\norm{x}\le\varkappa\l$, $\varkappa=\left(1-\dfrac{\varepsilon}{2}\right)[A(2k_2)^m]^{-1}$.
This gives the existence and unicity of 
$\zeta_0^*(Y;\l)\in I_x$ such that $X_0^*(Y;\l)=(0,\zeta_0^*(Y;\l))$ is the unique stationary
point of $F_2(X;Y;\l)$, with respect to $X$, which belongs to the support of $V_2(X;Y;\l)$
for $\norm{x}\le\varkappa\l$.

\noindent
First of all, the presence of the factors $H_1\!\left(\dfrac{\norm{x}\norm{(\l\zeta)^\frac{1}{m}\varsigma}^m}{\l}\right)$ and $H_2((\l\zeta)^\frac{1}{m})$ in the expression of $U_2$ imply
$(\l\zeta)^\frac{1}{m}\ge k_2\Rightarrow \norm{(\l\zeta)^\frac{1}{m}}\le(1+k_2^{-2})^\frac{1}{2}
(\l\zeta)^\frac{1}{m}$ and

\begin{align*}
	(2k_1)^{-1}\le&\,\frac{\norm{x}\norm{(\l\zeta)^\frac{1}{m}\varsigma}^m}{\l}
	\le \norm{x}(1+k_2^{-2})^\frac{m}{2}\zeta
	\\
	&\Rightarrow
	[2k_1(1+k_2^{-2})^\frac{m}{2}]^{-1}\le\norm{x}\zeta<
	\frac{\norm{x}\norm{(\l\zeta)^\frac{1}{m}\varsigma}^m}{\l}\le2k_1
	\\
	&\Rightarrow \zeta\in[c^{-1}\norm{x}^{-1},c\norm{x}^{-1}],
	\;c=2k_1(k_2^{-2}+1)^\frac{m}{2}.
\end{align*}

\noindent
Since $k_1>4AC>2AC$, clearly $I_x\subset[c^{-1}\norm{x}^{-1},c\norm{x}^{-1}]$. With
an arbitrarily chosen $\varepsilon\in\left(0,\dfrac{1}{2}\right)$, take
$k_2>\max\{B,1\}$ such that $\l\zeta>{k_2^m}$ implies $|S^{-1,0}(x,(\l\zeta)^\frac{1}{m})|\le\dfrac{\varepsilon}{2}$ and $|\zeta_0\zeta^{-1} S^{-1,0}(x,(\l\zeta)^\frac{1}{m})|\le k_0<1$, which is possible, in view of \eqref{disSG} and of the fact that $\zeta_0\zeta^{-1}$ is bounded
on $\supp V_2$. Fix $\l\ge\l_0>2k_1\norm{2k_2}^m$ and $\norm{x}\le\varkappa\l$.
Then, on $\supp V_2$, 
\begin{align*}
	\zeta\in I_x&\Rightarrow \l\zeta>
	\left(1-\dfrac{\varepsilon}{2}\right)^{-1}A(2k_2)^m
	\norm{x}\;
	A^{-1}\left(1-\dfrac{\varepsilon}{2}\right)\norm{x}^{-1}=(2k_2)^m>k_2^m
	\\
	&\Rightarrow G(\zeta;Y;\l)=\zeta_0(1+S^{-1,0}(x,(\l\zeta)^\frac{1}{m}))
	\in \left[
	A^{-1}\left(1-\dfrac{\varepsilon}{2}\right)\norm{x}^{-1},
	A\left(1+\dfrac{\varepsilon}{2}\right)\norm{x}^{-1}
	\right]= I_x
	\\
	&\Leftrightarrow G(.;Y;\l)\colon I_x \to I_x.
\end{align*}
Since $|\partial_\zeta G(\zeta;Y;\l)|=
|\zeta_0\zeta^{-1}S^{-1,0}(x,(\l\zeta)^\frac{1}{m})|\le k_0<1$, for all $\zeta\in I_x$,
$\norm{x}\le\varkappa\l$, we have proved that for any choice
of $Y\in \widetilde{S}_Y$, $\l\ge\l_0$ as above, $G(.;Y;\l)$ has a unique fixed point in
$\zeta_0^*=\zeta_0^*(Y;\l)\in I_x$, solution of $\zeta=G(\zeta;Y;\l)$.

By well-known corollaries of the Fixed Point Theorem for strict contractions
on compact subsets of metric spaces, we of course have that $\zeta_0^*$ depends
smoothly on $Y$ and $\l$. Moreover, since $\zeta_0^*\in I_x$
for all $Y\in \widetilde{S}_Y$, $\l\ge\l_0$, obviously $\zeta_0^*\sim\norm{x}^{-1}$ and
\[
	\zeta_0^*(Y;\l)=\zeta_0(1+S^{-1,0}(x,(\l\zeta_0^*(Y;\l))^\frac{1}{m})\to\zeta_0(\varsigma, x),
	\qquad\l\to+\infty,
\]
pointwise for any $(\varsigma, x)$. Moreover, by the choices of $k_1$, $k_2$ and $\varepsilon$,
\begin{align*}
	\frac{\norm{x}\norm{(\l\zeta_0^*(\varsigma, x;\l))^\frac{1}{m}}^m}{\l}
	&=\left[\frac{\norm{x}^\frac{2}{m}}{\l^\frac{2}{m}}+(\norm{x}\zeta_0^*)^\frac{2}{m}\right]^\frac{m}{2}
	>A^{-1}\left(1-\frac{\varepsilon}{2}\right)>k_1^{-1},
	\\
	\frac{\norm{x}\norm{(\l\zeta_0^*(\varsigma, x;\l))^\frac{1}{m}}^m}{\l}
	&<\left[\varkappa^\frac{2}{m}+ 
	\left(A\left(1+\frac{\varepsilon}{2}\right)\right)^\frac{2}{m}\right]^\frac{m}{2}
	=A\left[\left(1-\dfrac{\varepsilon}{2}\right)^\frac{2}{m}A^{-\frac{4}{m}}(2k_2)^{-2}
	+\left(1+\frac{\varepsilon}{2}\right)^\frac{2}{m}\right]^\frac{m}{2}
	\\
	&<k_1,
	\\
	\norm{x}\le\varkappa\l&\Leftrightarrow
	\l A^{-1}\left(1-\dfrac{\varepsilon}{2}\right)\norm{x}^{-1}\ge(2k_2)^m\Rightarrow
	\l \zeta_0^*(\varsigma, x;\l)\in [(2k_2)^m,+\infty).
\end{align*}
These imply, for any $\varsigma\in\Sn$, $x\in\R^n$, $\l\ge\l_0$
such that $\norm{x}\le\varkappa\l$,
\[
H_1\left(\dfrac{\norm{x}\norm{(\l\zeta_0^*(\varsigma, x;\l))^\frac{1}{m}}^m}{\l}\right)=1
\quad\mbox{ and }\quad
1-H_2((\l\zeta_0^*(\varsigma, x;\l))^\frac{1}{m})=1.
\]
Of course, by the choice of $H_3$, for $Y\in\widetilde{S}_Y$, $\l\ge\l_0$,
\[
\zeta_0^*\in I_x\Rightarrow H_3\left[\frac{\zeta_0^*(\varsigma, x;\l)}{\zeta_0(\varsigma,x)}-1\right]=1.
\]

                    %%%%%%%%%%%%%%%%%%%%%%%%%%%%%%%%%%%%%%%%%%%
                          %                Bibliografia             %
                            %%%%%%%%%%%%%%%%%%%%%%%%%%%%%%%%%%%%%%%%%%%
\addcontentsline{toc}{part}{References}

\bibliographystyle{abbrv}

\end{document}